\documentclass[11pt]{article} 
\usepackage{graphicx}
\pdfoutput=1
\usepackage{amssymb}
\usepackage{epstopdf}
\usepackage{color}
\usepackage{float}
\newcommand \RR {\mathbb R}

\newtheorem{property}{Property}

\long\def\symbolfootnote[#1]#2{\begingroup%
\def\thefootnote{\fnsymbol{footnote}}\footnote[#1]{#2}\endgroup}

\title{Convergent and conservative schemes \\
for nonclassical solutions based on kinetic relations. I.} 
\author{Benjamin Boutin$^{1,2}$, Christophe Chalons$^{1,3}$, 
Fr\'ed\'eric Lagouti\`ere$^{1,3}$, 
\\
and 
\\
Philippe G. LeFloch$^{1}$}

\date{} 

\begin{document}
\maketitle
\symbolfootnote[0]{$^1$ Laboratoire Jacques-Louis Lions \& Centre National de la
Recherche Scientifique, UMPC Universit\'e de Paris~6,
75252 Paris, France. 
}
\symbolfootnote[0]{$^2$ DEN/DANS/DM2S/SFME/LETR, CEA Saclay, 91191 Gif-sur-Yvette,
France.}
\symbolfootnote[0]{$^3$ Universit\'e Paris Diderot (Paris~7), 75251 Paris, France. 
\newline 
E-mail: Boutin@ann.jussieu.fr, Chalons@math.jussieu.fr, Lagoutie@math.jussieu.fr,
LeFloch@ann.jussieu.fr}

\begin{abstract}
We propose a new numerical approach to compute nonclassical solutions
to hyperbolic conservation laws. The class of finite difference
schemes presented here is fully conservative and keep nonclassical
shock waves as sharp interfaces, contrary to standard finite
difference schemes. The main challenge is to achieve, at the
discretization level, a consistency property with respect to a prescribed kinetic
relation. The latter is required for the selection of physically
meaningful nonclassical shocks. Our method is based on a
reconstruction technique performed in each computational cell that may
contain a nonclassical shock. To validate this approach, we establish
several consistency and stability properties, and we perform careful 
numerical experiments. The convergence of the algorithm toward the
physically meaningful solutions selected by a kinetic relation
is demonstrated numerically for several test cases, including
concave-convex as well as convex-concave flux-functions. 
\end{abstract}

\section{Introduction}
\label{0-0}

\subsection*{State of the art} 

We are interested here in the challenging issue of numerically
computing {\sl nonclassical} solutions (containing undercompressive
shocks) to nonlinear hyperbolic conservation laws. Nonclassical
solutions have the distinctive feature of being dynamically driven by
small-scale effects such as diffusion, dispersion, and other
high-order phenomena. Their selection requires an additional jump
relation, called a kinetic relation, and introduced in the context of
phase transition dynamics
\cite{sl1,sl2,tr1,tr2,ak1,ak2,fs,leflochARMA,hl1,hl2,shearer1,shearer2}, and 
investigated by LeFloch and collaborators in the
context of general hyperbolic systems of conservation laws (see \cite{lefloch1} for a review). 

From pioneering work by Hayes and LeFloch \cite{hl1,hl2} it is now
recognized that standard finite difference schemes do not converge to
nonclassical solutions selected by the prescribed kinetic function. In
fact, kinetic functions can be associated not only with continuous
models, but with the finite difference schemes themselves. Achieving a
good agreement between the continuous and the numerical kinetic
functions has been found to be very challenging. 

In the present paper, we will show how to enforce the validity of the
kinetic relation at the numerical level, and we design a {\sl fully
conservative} scheme which combines the advantages of standard finite
differences and Glimm-type (see below) approaches. 

Nonclassical shocks and other
phase transitions are naturally present in many models of continuum
physics, especially in the modeling of real fluids governed by complex
equations of state. This is the case, for instance, of models
describing the dynamics of liquid-vapor phase transitions in compressible
fluids, or of solid-solid phase transformations in materials such as
memory alloys. For numerical work in this direction we refer to
\cite{hlr1,hlz1,cl2,mr2,mr1}. 


\subsection*{Setting for this paper} 

We restrict here attention to scalar conservation laws  
\begin{equation} 
\label{1}
\begin{array}{l}
{\partial}_t u + {\partial}_x f(u) = 0, \qquad u(x,t) \in \mathbb{R},
\,\,\, (x,t) \in \mathbb{R} \times \mathbb{R}^{+},\\
u(x,0) = u_0(x), 
\end{array}
\end{equation}
and postpone the discussion of systems of conservation laws to the
follow-up paper \cite{BCLL-two}. The above equation must be
supplemented with an entropy inequality of the form 
\begin{equation} 
\label{2}
{\partial}_t U(u) + {\partial}_x F(u) \leq 0. 
\end{equation}
Here, $t$ denotes the time variable, $x$ the (one-dimensional) space
variable, $f: \mathbb{R} \to \mathbb{R}$ the flux function, and
$(U,F)$ is any strictly convex mathematical entropy pair. That is, $U:
\mathbb{R} \to \mathbb{R}$ is strictly convex and $F: \mathbb{R} \to
\mathbb{R}$ is given by $F'=U'f'$. Equations (\ref{1}) and (\ref{2})
are imposed in the distributional sense.

We rely here on the theory of nonclassical solutions based on kinetic relations, 
established in \cite{lefloch1}. The flux $f$ is assumed to be {\sl nonconvex,} which is the
source of mathematical and numerical difficulties. From the
mathematical standpoint, a single entropy inequality like (\ref{2})
does not suffice to select a unique solution. This can be seen already
at the level of the Riemann problem, corresponding to
(\ref{1})-(\ref{2}) when $u_0$ has the piecewise constant form 
\begin{equation} \label{ditr}
u_0(x)
=
\left\{
\begin{array}{rcl}
u_l, & & x < 0, \\
u_r, & & x > 0, \\
\end{array}
\right. 
\end{equation}
$u_l$ and $u_r$ being constant states. The Riemann problem admits (up
to) a one-parameter family of solutions (see Chapter 2 in
\cite{lefloch1}). However, these solutions contain discontinuities
violating the standard Lax shock inequalities, which are referred to
as {\sl nonclassical}. They are essential from the physical
standpoint, and should be retained. This non-uniqueness can be
fixed however, provided an additional algebraic condition, the
so-called {\sl kinetic relation,} is imposed on each nonclassical
shock. Consider a shock connecting a left-hand state $u_-$ to a
right-hand state $u_+$ and propagating with the speed $\sigma$ given
by the usual Rankine-Hugoniot relation, that is, 
\begin{equation} \label{ss}
u(x,t) = \left\{
\begin{array}{rcl}
u_-, &  & x < \sigma t, \\
u_+, & & x > \sigma t, \\
\end{array}
\right. 
\quad   \quad 
\sigma = \sigma(u_-,u_+)=\frac{f(u_+)-f(u_-)}{u_+-u_-}.
\end{equation} 
The kinetic relation takes the form 
\begin{equation} 
\label{kc}
u_+ = \varphi^{\flat}(u_-) \,\,\,\, \mbox{for all nonclassical shocks},
\end{equation}
where $\varphi^{\flat}$ is the so-called {\sl kinetic
function}. Equivalently, denoting by $\varphi^{-\flat}$ the inverse of
the kinetic function it may be preferable to write $u_- =
\varphi^{-\flat}(u_+)$. The kinetic relation implies that the
right-hand (respectively left-hand) state is no longer free (as in a
classical shock wave) but depends explicitly on the left-hand
(respectively right-hand) state. 


\subsection*{Objectives in this paper} 

At the numerical level, several strategies exist in the literature in
order to take into account the kinetic relation (\ref{kc}). We can
distinguish between diffuse interface methods and sharp interface
methods. 

In the first approach, one assumes that the kinetic relation is
derived from an augmented continuous model and, in order to take into
account the internal structure of nonclassical discontinuities, one 
attempts to resolve the effects dues to (small) diffusive and
dispersive terms that generate them. It is then possible to construct
conservative schemes that mimic at the numerical level the effect of
the regularized models. Due to the great sensitivity of nonclassical
solutions with respect to small scales and numerical diffusion, it
turns out that numerical results are satisfactory for shocks with
moderate amplitude, but discrepancies between the exact and the
numerical kinetic function arise with shocks with large amplitudes and
in long-time computations. For this circle of ideas we refer the
reader to \cite{hl1,hl2}, and the follow-up papers \cite{lr1,cl1,cl2}.  

In the second approach, small scale features are not explicitly taken
into account. Instead, the kinetic relation is included, in a
way or another, in the design of the numerical scheme. This is the
case of the random choice and front tracking schemes. It should be
mentioned here that the Glimm scheme and front tracking schemes do
converge to exact solutions even in presence of nonconclassical
shocks; see \cite{leflochARMA,lefloch1,ls1} for the theoretical
aspects and Chalons and LeFloch \cite{cl4} for a numerical study of
the Glimm scheme. These schemes require the explicit knowledge of the
underlying nonclassical Riemann solver, which may be expensive
numerically, and this motivated the introduction of the 
so-called transport-equilibrium scheme by Chalons \cite{chalons2,chalons3}. 

In \cite{hlz1}, Hou, LeFloch, and Zhong proposed a class of converging schemes  
for the computation of propagating solid-solid phase boundaries. 
More recently, Merckle and Rohde \cite{mr1}
developed a ghost-fluid type algorithm for a model of dynamics of
phase transition. These schemes provide satisfactory numerical results, as
nonclassical discontinuities are sharply and accurately
computed. Although the convergence of the methods was demonstrated
numerically, their main drawback in practice is similar to the
Glimm-type schemes and the property of strict conservation of the
conservative variable $u$ fails. 

Building on these previous works, our objective in this paper is to
design a fully conservative, finite difference scheme for the
approximation of nonclassical solutions to the hyperbolic conservation
law (\ref{1}). Our basic strategy relies on the discontinuous
reconstruction technique proposed recently in Lagouti\`ere
\cite{lagoutiere3,lagoutiere4} which has been found to be particularly
efficient to computing {\sl classical} solutions of (\ref{1}) with
moderate numerical diffusion. 

In our approach below, the kinetic function $\varphi^{\flat}$ is
included explicitly in the algorithm, in such a way that nonclassical shocks are computed
(essentially) {\sl exactly} while classical shocks suffer moderate 
numerical diffusion. To validate our strategy we perform various 
numerical experiments and, in particular, draw the kinetic function
associated with our scheme. As the mesh is refined, we observe that
the approximate kinetic function converges toward the analytic kinetic
function. The scheme also enjoys several fundamental stability
properties of consistency with the conservative form of the
equation and (like the Glimm scheme) with single nonclassical discontinuities. 


\section{Nonclassical Riemann solver with kinetics} 
 \label{2-0}

\subsection*{Assumption on the flux-function}

We describe here the nonclassical Riemann solver introduced and
investigated in LeFloch \cite{lefloch1}. Note in passing that this
solver was later extended in \cite{ls1} to include also a nucleation
criterion. 

Consider the problem (\ref{1})-(\ref{2})-(\ref{kc}) for a given
Riemann initial data (\ref{ditr}). Throughout this paper we assume
that the flux $f$ is either {\sl concave-convex} or {\sl
convex-concave}, that is, satisfies the conditions (for all $u \neq
0$) 
\begin{equation} 
\label{conca-conve}
\begin{array}{c}
u f''(u) > 0, \qquad f'''(0) \neq 0,
\qquad 
\lim_{|u| \to +\infty} f'(u) = + \infty,
\end{array}
\end{equation}
or
\begin{equation} \label{conve-conca}
\begin{array}{c}
u f''(u) < 0, \qquad f'''(0) \neq 0,
\qquad
\lim_{|u| \to +\infty} f'(u) = - \infty, 
\end{array}
\end{equation}
respectively. The functions $f(u) = u^3 + u$ and $f(u) = -u^3 - u$ 
are prototypes of particular interest, used later in this paper 
for the validation of the proposed numerical strategy. 

Let $\varphi^{{\natural}} : \mathbb{R} \to \mathbb{R}$ be the unique
function defined by $\varphi^{{\natural}}(0)=0$ and for all $u \neq
0$, $\varphi^{{\natural}}(u) \neq u$ is such that the line passing
through the points $(u,f(u))$ and
$(\varphi^{{\natural}}(u),f(\varphi^{{\natural}}(u)))$ is tangent to
the graph of $f$ at point 
$(\varphi^{{\natural}}(u),f(\varphi^{{\natural}}(u)))$: 
$$
f'(\varphi^{{\natural}}(u)) =
\frac{f(u)-f(\varphi^{{\natural}}(u))}{u-\varphi^{{\natural}}(u)}. 
$$
This function is smooth, monotone decreasing and onto thanks to
(\ref{conca-conve}) or (\ref{conve-conca}). We denote by
$\varphi^{-{\natural}} : \mathbb{R} \to \mathbb{R}$ its inverse
function. \\


\subsection*{Concave-convex flux functions}

Let us assume that $f$ obeys (\ref{conca-conve}) and let
$\varphi^{\flat} : \mathbb{R} \to \mathbb{R}$ be a kinetic function,
that is (by definition) a monotone decreasing and Lipschitz continuous
mapping such that 
\begin{equation} \label{hypkf1}
\begin{array}{rcl}
\varphi^{\flat}_0(u) < \varphi^{\flat}(u) \leq \varphi^{{\natural}}(u),
&  & u > 0, \\ 
\varphi^{{\natural}}(u) \leq \varphi^{\flat}(u) < \varphi^{\flat}_0(u),
&  & u < 0. 
\end{array}
\end{equation}
From $\varphi^{\flat}$, we define the function 
$\varphi^{\sharp} : \mathbb{R} \to \mathbb{R}$ such that the line
passing through the points $(u,f(u))$ and
$(\varphi^{\flat}(u),f(\varphi^{\flat}(u)))$ with $u \neq 0$ also cuts
the graph of the flux function $f$ at point
$(\varphi^{\sharp}(u),f(\varphi^{\sharp}(u)))$ with
$\varphi^{\sharp}(u) \neq u$ and $\varphi^{\sharp}(u) \neq
\varphi^{\flat}(u)$: 
$$
\frac{f(u)-f(\varphi^{\flat}(u))}{u-\varphi^{\flat}(u)} 
= 
\frac{f(u)-f(\varphi^{{\sharp}}(u))}{u-\varphi^{{\sharp}}(u)}.
$$
The nonclassical Riemann solver associated with
(\ref{1})-(\ref{2})-(\ref{ditr})-(\ref{kc}) is given as follows. \\ 
When $u_l > 0$: 

$(1)$ If $u_r \geq u_l$, the solution is a rarefaction wave connecting 
$u_l$ to $u_r$.

$(2)$ If $u_r \in [{\varphi}^{\sharp}(u_l), u_l)$, the solution is a
classical shock wave connecting $u_l$ to $u_r$.

$(3)$ If $u_r \in ({\varphi}^{\flat}(u_l), {\varphi}^{\sharp}(u_l))$,
the solution contains a nonclassical shock connecting $u_l$ to
${\varphi}^{\flat}(u_l)$, followed by a classical shock connecting
${\varphi}^{\flat}(u_l)$ to $u_r$. 

$(4)$ If $u_r \leq {\varphi}^{\flat}(u_l)$, the solution contains a
nonclassical shock connecting $u_l$ to ${\varphi}^{\flat}(u_l)$,
followed by a rarefaction connecting ${\varphi}^{\flat}(u_l)$ to
$u_r$. \\ 
When $u_l \leq 0$: 

$(1)$ If $u_r \leq u_l$, the solution is a rarefaction wave connecting
$u_l$ to $u_r$.

$(2)$ If $u_r \in [u_l, {\varphi}^{\sharp}(u_l))$, the solution is a
classical shock wave connecting $u_l$ to $u_r$.

$(3)$ If $u_r \in ({\varphi}^{\sharp}(u_l), {\varphi}^{\flat}(u_l))$,
the solution contains a nonclassical shock connecting $u_l$ to
${\varphi}^{\flat}(u_l)$, followed by a classical shock connecting
${\varphi}^{\flat}(u_l)$ to $u_r$. 

$(4)$ If $u_r \geq {\varphi}^{\flat}(u_l)$, the solution contains a
nonclassical shock connecting $u_l$ to ${\varphi}^{\flat}(u_l)$,
followed by a rarefaction connecting ${\varphi}^{\flat}(u_l)$ to
$u_r$. 


\subsection*{Convex-concave flux functions}

We next assume that $f$ satisfies the condition
(\ref{conve-conca}). Let $\varphi^{\flat} : \mathbb{R} \to \mathbb{R}$
be a kinetic function, that is, a monotone decreasing and Lipschitz
continuous map such that 
\begin{equation} \label{hypkf2}
\begin{array}{rcl}
\varphi^{\flat}_0(u) < \varphi^{\flat}(u) \leq
\varphi^{-{\natural}}(u), & & u < 0, \\ 
\varphi^{-{\natural}}(u) \leq \varphi^{\flat}(u) <
\varphi^{\flat}_0(u), &  & u > 0. 
\end{array}
\end{equation}
We then define $\rho(u,v) \in \mathbb{R}$ if 
$v \neq u$ and $v \neq \varphi^{{\natural}}(u)$ by 
$$
\frac{f(\rho(u,v))-f(u)}{\rho(u,v) - u} 
= 
\frac{f(v)-f(u)}{v-u}
$$
with $\rho(u,v) \neq u$ and $\rho(u,v) \neq v$, and extend the 
function $\rho$ by continuity otherwise. Note that
$\varphi^{{\sharp}}(u) = \rho(u,\varphi^{\flat}(u))$ where 
$\varphi^{{\sharp}}$ is defined as in the case of a concave-convex
flux function. The nonclassical Riemann solver associated with
(\ref{1})-(\ref{2})-(\ref{ditr})-(\ref{kc}) is given as follows. \\ 
When $u_l > 0$: 

$(1)$ If $u_r \geq u_l$, the solution is a classical shock connecting 
$u_l$ to $u_r$.

$(2)$ If $u_r \in [0, u_l)$, the solution is a rarefaction wave 
connecting $u_l$ to $u_r$.

$(3)$ If $u_r \in ({\varphi}^{\flat}(u_l), 0)$, the solution contains
a rarefaction wave connecting $u_l$ to ${\varphi}^{-\flat}(u_r)$,
followed by a nonclassical shock connecting ${\varphi}^{-\flat}(u_r)$
to $u_r$. 

$(4)$ If $u_r \leq {\varphi}^{\flat}(u_l)$, the solution contains: 

$\quad$ (i) a classical shock connecting $u_l$ to
${\varphi}^{-\flat}(u_r)$, followed by a nonclassical shock connecting
${\varphi}^{-\flat}(u_r)$ to $u_r$, if $u_l >
\rho({\varphi}^{-\flat}(u_r),u_r)$. 

$\quad$ (ii) a 
classical shock connecting $u_l$ to $u_r$, if
$u_l \leq \rho({\varphi}^{-\flat}(u_r),u_r)$. \\
When $u_l \leq 0$:

$(1)$ If $u_r \leq u_l$, the solution is a classical shock connecting 
$u_l$ to $u_r$.

$(2)$ If $u_r \in (u_l, 0]$, the solution is a rarefaction wave 
connecting $u_l$ to $u_r$.

$(3)$ If $u_r \in (0, {\varphi}^{\flat}(u_l))$, the solution contains
a rarefaction wave connecting $u_l$ to ${\varphi}^{-\flat}(u_r)$,
followed by a non classical shock connecting ${\varphi}^{-\flat}(u_r)$
to $u_r$. 

$(4)$ If $u_r \geq {\varphi}^{\flat}(u_l)$, the solution contains: 

$\quad$ (i) a 
classical shock connecting $u_l$ to ${\varphi}^{-\flat}(u_r)$,
followed by a nonclassical shock connecting ${\varphi}^{-\flat}(u_r)$
to $u_r$, if $u_l < \rho({\varphi}^{-\flat}(u_r),u_r)$. 

$\quad$ (ii) a classical shock connecting $u_l$ to $u_r$, if
$u_l \geq \rho({\varphi}^{-\flat}(u_r),u_r)$.

Observe that the convex-concave case can in principle be deduced from
the concave-convex case, by replacing $f$ by $-f$ and $x$ by
$-x$. Nevertheless, it is useful to keep the above two descriptions
in mind, since there is a dramatic difference between the Riemann
solvers: the nonclassical shock always connects $u_l$ to
${\varphi}^{\flat}(u_l)$ in the concave-convex case, and
${\varphi}^{-\flat}(u_r)$ to $u_r$ in the convex-concave case. The
numerical method we are going to describe must take this feature into
account, and as we will explain it is necessary to take into account
both ${\varphi}^{\flat}$ and ${\varphi}^{-\flat}$ in the design of the
scheme. 


\section{Motivations and difficulties} 
\label{5-0}

\subsection*{Notation} 

Our aim is to design a scheme for the numerical approximation of the
nonclassical solutions to (\ref{1})-(\ref{2})-(\ref{kc}). To this end,
we consider the general class of finite volume methods. Introducing
constant space and time lengths $\Delta x$ and $\Delta t$ for the
space and time discretization, we can set $x_{j+1/2} = j \Delta x$, $j
\in \mathbb{Z}$, and $t^n = n \Delta t$, $n \in \mathbb{N}$. The
discretization consists, at each time $t^n$, of a piecewise constant
function $x \mapsto {u}_{\nu}(x,t^n)$ which should be an approximation
of the exact solution ${u}(x,t^n)$ on the cell
$\mathcal{C}_{j}=[x_{j-1/2};x_{j+1/2})$: 
$$
{u}_{\nu}(x,t^n) = {u}^n_{j}, \qquad 
x \in C_{j}, \,\,\, j \in \mathbb{Z}, \,\,\, n \in
\mathbb{N}.
$$
Here, $\nu$ refers to the ratio $\Delta t / \Delta x$.
The initial data at the time $t=0$ is denoted by $u_0$ and we define
the sequence $({u}^0_{j})_{j \in \mathbb{Z}}$: 
\begin{equation} 
\label{initialisation}
{u}^0_{j}= \frac{1}{\Delta x} \int_{x_{j-1/2}}^{x_{j+1/2}} {u}_0(x)
dx, \qquad 
j \in \mathbb{Z}. 
\end{equation} 

The starting point in the conception of our algorithm is a 
few conventional interpretation of the constant values 
${u}^n_j$, $j \in \mathbb{Z}$. As suggested by the proposed
initialization (\ref{initialisation}), ${u}^n_j$ is usually, and
rightly, seen as an approximate value of the average 
on cell $\mathcal{C}_{j}$ 
of the exact solution at time $t^n$. Integrating equation (\ref{1})
over the slab $\mathcal{C}_j \times [t^n, t^{n+1}]$ and using Green's
formula, it is thus natural to define $({u}^{n+1}_j)_j$ from
$({u}^n_j)_j$ and a conservative scheme of the following form
\begin{equation} \label{schema_cons}
{u}^{n+1}_j = {u}^{n}_j - \frac{\Delta t}{\Delta x} 
(f^n_{j+1/2} - f^n_{j-1/2}), \,\,\, \,\,\,j \in
\mathbb{Z}, 
\end{equation}
where $f^n_{j+1/2}$ represents an approximate value of the flux that
passes through the interface $x_{j+1/2}$ between the times $t^n$ and
$t^{n+1}$. \\ 
Here, we shall also consider ${u}^n_j$ as a given information, on cell
$\mathcal{C}_j$ and at time $t^n$, on the structure of the exact
Riemann solution associated with inital states $u_l = {u}^n_{j-1}$ and
$u_r = {u}^n_{j+1}$ which will develop at the next times $t > t^n$. At
this stage, one easily realize that if this information is precise ({\sl i.e.} close to what will really happen), then we should be 
be in a good position to define accurately the numerical fluxes $f^n_{j+1/2}$ and then predict the approximate values of the solution
at time $t^{n+1}$.


\subsection*{Linear advection equation}

\begin{figure}[htbp]
\begin{center}
\includegraphics[width=8.cm]{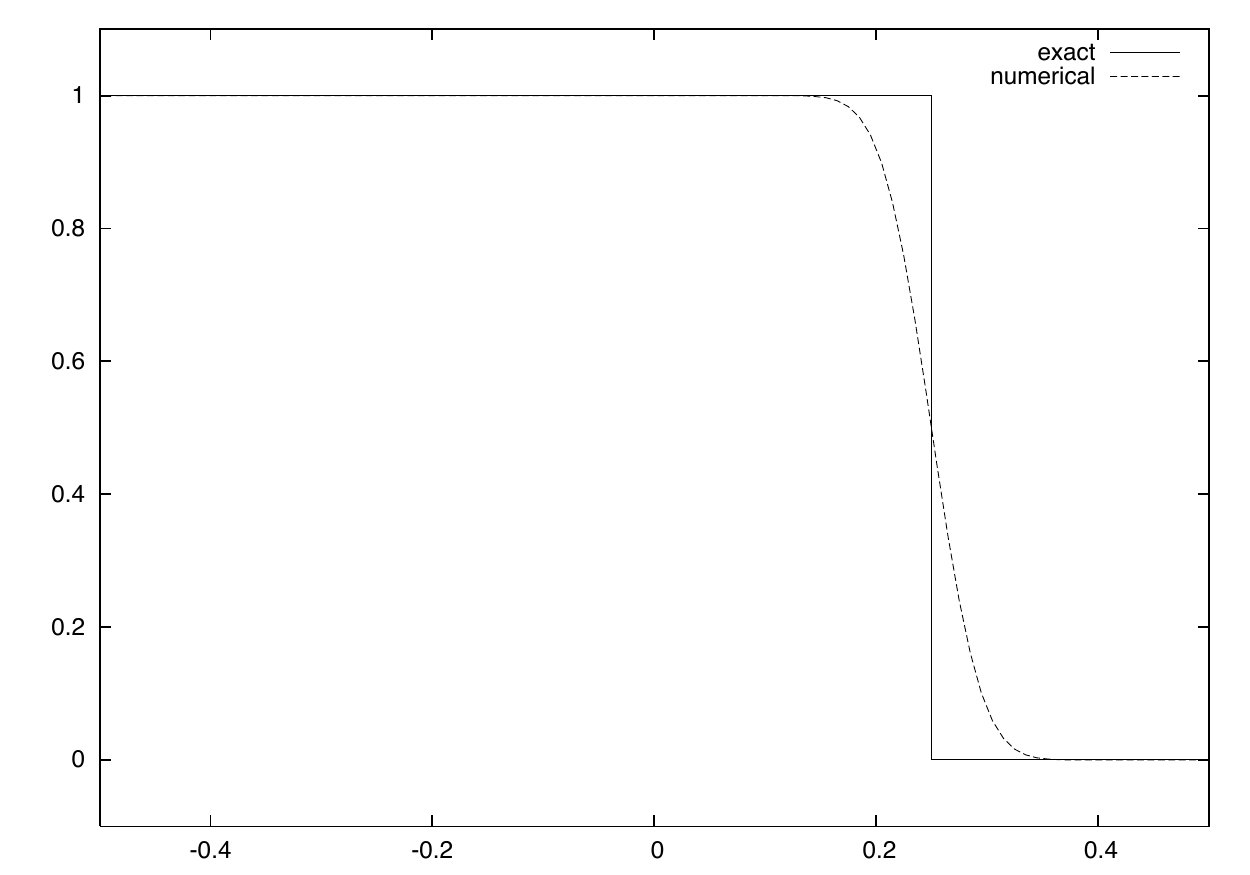}
\includegraphics[width=8.cm]{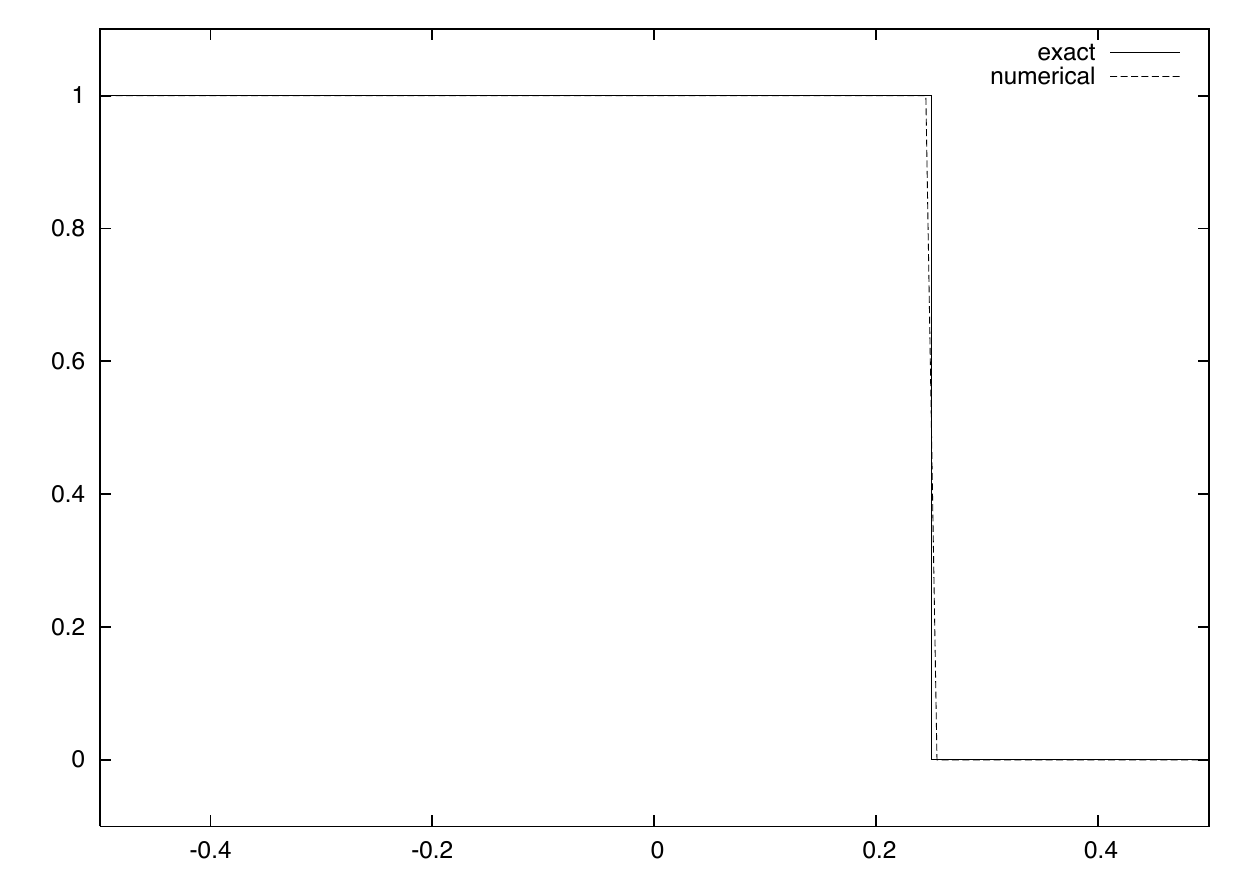}
\vspace{-0.45cm}
\caption{Linear advection - upwind scheme (left-hand) and reconstruction
scheme (right-hand).} 
\label{advection12}
\end{center}
\end{figure}

As a first illustration, let us consider the linear advection with
constant velocity $a > 0$, that is, the scalar conservation law with
flux $f(u) = a u$. In this case, the weak solution to the
initial-value problem for (\ref{1}) is unique, and is given explicitly
as $u(t,x) = u_0(x-at)$. Hence, neither the entropy condition
(\ref{2}) nor the kinetic condition (\ref{kc}) are necessary. The
basic scheme for approximating this solution is the so-called upwind
scheme and corresponds to the choice $f^n_{j+1/2} = a u^n_j$ for all
$j \in \mathbb{Z}$. Recall that the CFL condition $a \Delta t / \Delta
x \leq \alpha$ for a given $\alpha \leq 1$ is mandatory for the
stability of the procedure. Figure \ref{advection12} (left-hand) shows the
corresponding numerical solution at time $t=0.25$ for $a = 1$, $\alpha
= 0.5$ and $u_l = 1$, $u_r = 0$ in (\ref{ditr}). The mesh contains
$100$ points per unit. 
We observe that the numerical solution presents a 
good agreement with the exact one but contains numerical diffusion.
We propose the following interpretation. In some sense, the value
$u^n_j$ that we consider as an information on the Riemann solution
associated with initial states $u_l = u^n_{j-1}$ and $u_r = u^n_{j+1}$
is sufficient to correctly approach this solution when defining
$f^n_{j+1/2} = a u^n_j$, but not enough to avoid the numerical
diffusion. Note that the latter is expected but not hoped. In the
present situation, the fact is that we actually know what will happen
in the future, namely a propagation of the Riemann initial states
($u_l = u^n_{j-1}$ and $u_r = u^n_{j+1}$) with speed $a$. In
particular, no value different from $u^n_{j-1}$ and $u^n_{j+1}$ is
created so that information given by $u^n_j$ is clearly not
optimal. In the process of calculation of the numerical flux
$f^n_{j+1/2}$, we are thus tempted to add more information in the cell 
$\mathcal{C}_j$ when replacing, as soon as possible,
the constant state $u^n_j$ with a
discontinuity separating $u^n_{j-1}$ on the left and $u^n_{j+1}$ on
the right, and located at point $\overline{x}_j \in \mathcal{C}_j$. In
the forthcoming developments, the left and right states of this
reconstructed discontinuity will be noted $u^n_{j,l}$ and $u^n_{j,r}$,
respectively. Hence, we have here 
\begin{equation} \label{first_reconst}
u^n_{j,l} = u^n_{j-1}, 
\quad \quad
u^n_{j,r} = u^n_{j+1}.
\end{equation}  
We claim that this provides better information for calculating 
$f^n_{j+1/2}$ than the original one. Such a reconstruction is due to
conserve $u$ in order to be relevant, which defines $\overline{x}_j$
by the following constraint 
$$
( \overline{x}_j - x_{j-1/2} ) u^n_{j,l} + 
( x_{j+1/2} - \overline{x}_j ) u^n_{j,r} = 
( x_{j+1/2} - x_{j-1/2} ) u^n_j 
$$
which equivalently recast as
\begin{equation} \label{x_disc}
\overline{x}_j = x_{j-1/2} + 
\frac{u^n_{j,r} - u^n_{j}}{u^n_{j,r} - u^n_{j,l}} \Delta x.
\end{equation}
Then, the reconstruction is possible provided we have $0 \leq d^n_j
\leq 1$, with 
\begin{equation} \label{defdnj}
d^n_j = \frac{u^n_{j,r} - u^n_{j}}{u^n_{j,r} - u^n_{j,l}}.
\end{equation}

\noindent Now, let us introduce $\Delta t_{j+1/2}$ the time needed by
the reconstructed discontinuity to reach the interface $x_{j+1/2}$
(recall that $a >0$). We clearly have 
$$
\Delta t_{j+1/2} = \frac{1 - d^n_j}{a} \Delta x.
$$ 
In this case, the flux that passes through $x_{j+1/2}$ between times
$t^n$ and $t^{n+1} = t^n + \Delta t$ equals $f(u^n_{j,r})$ until $t^n
+ \Delta t_{j+1/2}$, and $f(u^n_{j,l})$ after (if $\Delta t_{j+1/2} <
\Delta t$). Therefore, we propose to set now
$$
\Delta t f^n_{j+1/2} = 
\min(\Delta t_{j+1/2}, \Delta t) f(u^n_{j,r}) +
\max(\Delta t - \Delta t_{j+1/2}, 0) f(u^n_{j,l}).  
$$ 
On Figure \ref{advection12} (right-hand), we have plotted the 
numerical solution given by this new numerical flux, leading to the 
so-called reconstruction scheme.
The parameters of the simulation are the same than those of
Figure \ref{advection12} (left-hand). 
We see that the more precise informations we have
brought on each cell $\mathcal{C}_j$ for calculating the numerical
fluxes make the scheme less diffusive than the original one. This
strategy was proposed (and is discussed in further details) in
\cite{lagoutiere3,lagoutiere4} (see also \cite{dl1,lagoutiere1}). 
In particular, it is shown therein that the numerical solution presented in Figure \ref{advection12} 
(right-hand) is {\sl exact} in the sense that $u^n_j$ equals the average of
the exact solution on $\mathcal{C}_j$, that is 
\begin{equation} 
\label{se}
u^n_j = \frac{1}{\Delta x} \int_{x_{j-1/2}}^{x_{j+1/2}} u(x,t^n) dx,
\qquad \, j \, \in \, \mathbb{Z}, \,\, \, n \, \in \, \mathbb{N}. 
\end{equation}
The corresponding numerical discontinuity separating $u_l$ and $u_r$
in then diffused on one cell at most. 


\subsection*{Godunov scheme with a nonclassical Riemann solver}

\begin{figure}[htbp]
\begin{center}
\includegraphics[width=8.0cm]{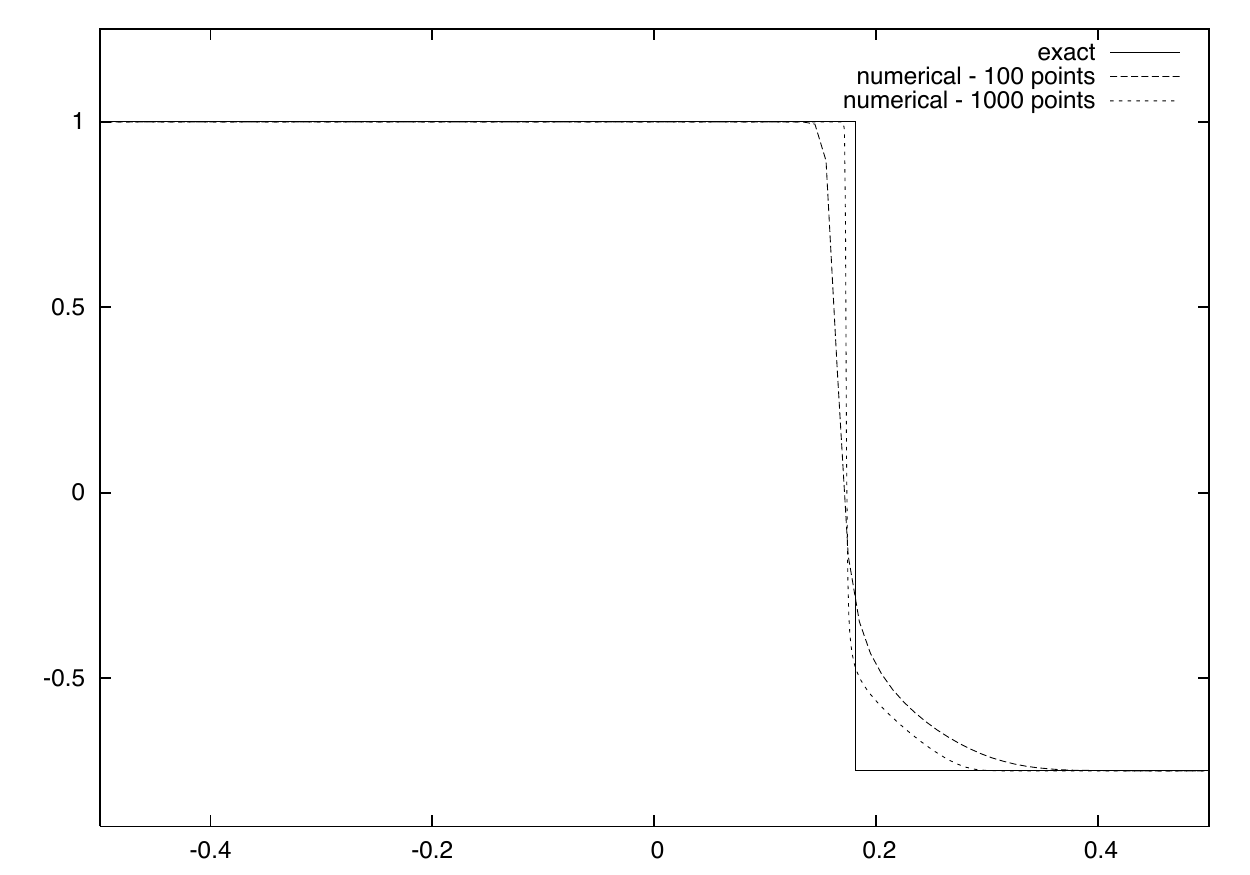}
\includegraphics[width=8.0cm]{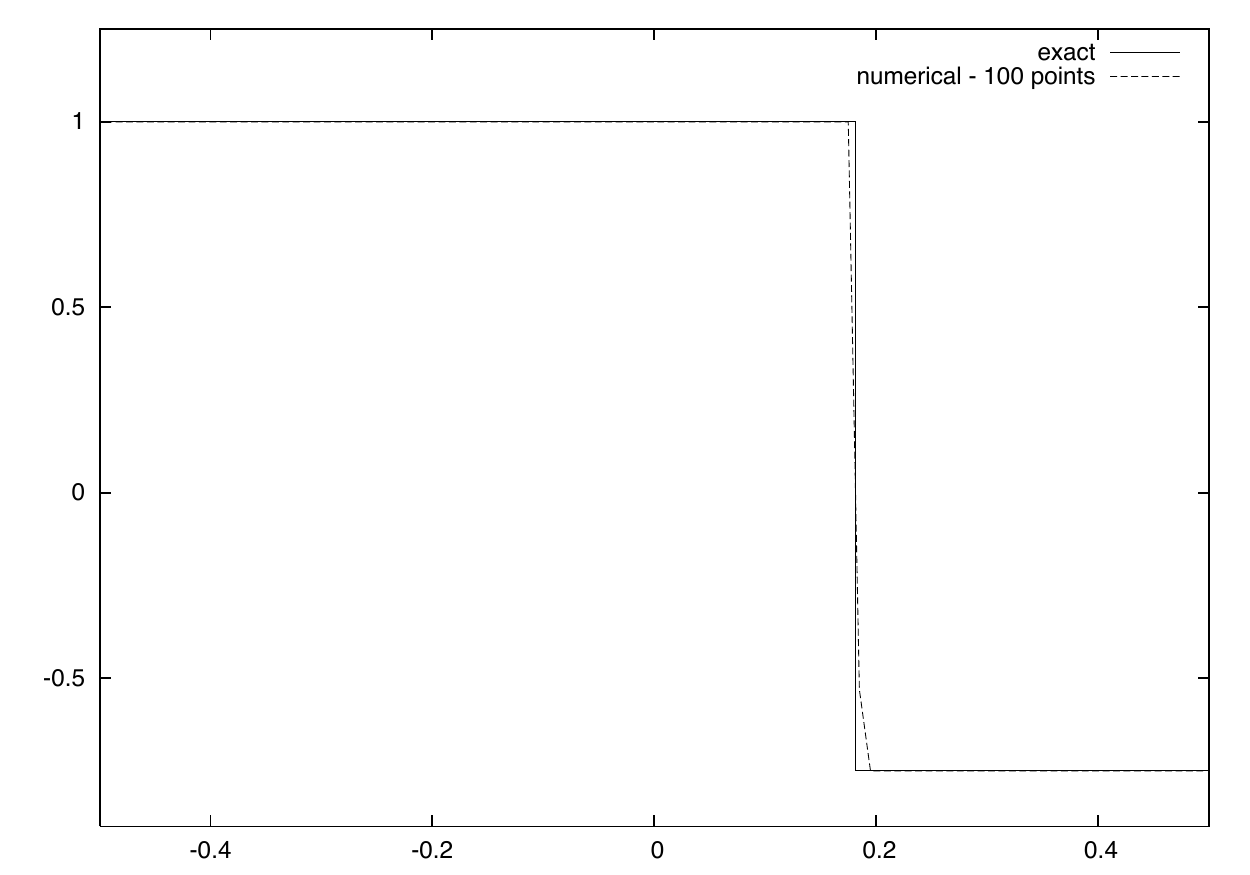}
\vspace{-0.45cm}
\caption{Propagating nonclassical shock - upwind scheme (left-hand) and
reconstruction scheme (right-hand).} 
\label{advection34}
\end{center}
\end{figure}

As a second illustration, let us go back to the problem 
(\ref{1})-(\ref{2})-(\ref{kc}) with a general 
concave-convex (or convex-concave) flux function $f$ with however, for
the sake of clarity, 
\begin{equation} \label{assumpfp}
f'(u) \geq 0, \qquad u \in \RR.
\end{equation} 
Here, we focus ourselves on a particular Riemann initial data 
(\ref{ditr}) such that $u_r = \varphi^{\flat}(u_l)$. In other words, 
the kinetic criterion is imposed on the initial discontinuity. The 
exact solution then corresponds to the propagation of this
discontinuity with speed $\sigma(u_l,u_r) > 0$ given by
Rankine-Hugoniot relation: 
\begin{equation} \label{srh}
 \sigma(u_l,u_r) = \frac{f(u_r) - f(u_l)}{u_r - u_l}.
\end{equation}
Figure \ref{advection34} (left-hand) represents the numerical solution
given by the upwind scheme $f^n_{j+1/2} = f(u^n_j)$ at time $t=0.1$,
for $f(u) = u^3 + u$ and $u_l = 1$. The kinetic function is taken to
be $\varphi^{\flat}(u) = - 0.75 \, u$ so that $u_r = - 0.75$. 

We observe a strong disagreement between the numerical solution and
the exact one. Indeed, the former is made of a (classical) shock
followed by a rarefaction wave while the latter is a single
(nonclassical) shock from $u_l$ to $u_r$. It is then clear that the
usual upwind scheme (as many others actually) is not adapted for the
computation of nonclassical solutions. The next result states that the
upwind scheme always converges towards the classical solution of
(\ref{1})-(\ref{2}). This scheme is then adapted for the computation
of classical solutions only. 

\

{\bf Property.}
{\sl 
Assume that $u_0 \in L^{\infty}(\mathbb{R})$ and $f$ is a smooth
function satisfying (\ref{assumpfp}). Then, under the CFL condition 
$$
\frac{\Delta t}{\Delta x} \, \max_{u \in A}
|f'(u)| \leq 1, 
$$
with $A := [\min_x u_0(x), \max_x u_0(x)]$
the upwind conservative scheme (\ref{schema_cons}) with $f^n_{j+1/2} =
f(u^n_j)$ converges towards the unique classical solution of
(\ref{1})-(\ref{2}). 
}

\ 

To establish this property, we only need to observe that, under the
assumption (\ref{assumpfp}) (propagation is only in one direction),
the upwind scheme is equivalent to the standard Godunov scheme
associated with the {\sl classical} Riemann solver of
(\ref{1})-(\ref{2})  
Then, standard compactness and consistency arguments apply and allow
us to conclude that the scheme converges towards the unique classical
solution. 

Obviously, the above property also holds if $f$ is assumed to be
decreasing if we define $f^n_{j+1/2} = f(u^n_{j+1})$. 


\section{A conservative scheme for nonclassical entropy solutions}

\subsection*{Preliminaries} 
In view of the discussion in the previous section and in order to
better evaluate the numerical fluxes $f^n_{j+1/2}$, let us obtain
some information beyond $u^n_j$ on cell $\mathcal{C}_j$. In the 
present instance of an isolated propagating discontinuity, it is
expected that the Riemann solution associated with initial states
$u^n_{j-1}$ and $u^n_{j+1}$ simply propagates the initial
discontinuity. This is actually true if $u^n_{j-1} = u_l$ and
$u^n_{j+1} = \varphi^{\flat}(u_l)$, or more generally if $u^n_{j+1} =
\varphi^{\flat}(u^n_{j-1})$. So that here again, we propose to replace
the constant state $u^n_j$ with a discontinuity separating $u^n_{j,l}$
and $u^n_{j,r}$ and located at point $\overline{x}_j$ given by
(\ref{x_disc}), as soon as possible {\sl i.e.} when $0 \leq d^n_j \leq
1$. We take 
\begin{equation} \label{second_reconst}
u^n_{j,l} = \varphi^{-\flat}(u^n_{j+1})
\quad \mbox{and} \quad 
u^n_{j,r} = \varphi^{\flat}(u^n_{j-1}).
\end{equation}
Note that this reconstruction is equivalent to (\ref{first_reconst}) 
provided that 
$u^n_{j-1} = u_l$ and $u^n_{j+1} = \varphi^{\flat}(u_l)$, or more
generally $u^n_{j+1} = \varphi^{\flat}(u^n_{j-1})$. Then, under the
assumption (\ref{assumpfp}), we again naturally set
$$
\Delta t f^n_{j+1/2} = 
\min(\Delta t_{j+1/2}, \Delta t) f(u^n_{j,r}) +
\max(\Delta t - \Delta t_{j+1/2}, 0) f(u^n_{j,l})
$$
with now
\begin{equation} \label{deltat_int_2}
\Delta t_{j+1/2} = \frac{1 - d^n_j}{\sigma(u^n_{j,l},u^n_{j,r})}
\Delta x. 
\end{equation} 
Figure \ref{advection34} (right-hand) highlights the benefit of such a
reconstruction. The numerical solution now fully agrees with the exact
one and is moreover free of numerical diffusion (the profile is
composed of a single point). We will show below that it is exact in
this case, in the sense that (\ref{se}) is still valid as in the
linear case. 


\subsection*{The scheme} 

On the basis of the above motivations and illustrations, 
we follow the description of our algorithm by
considering the general situation. Assuming as given a sequence
$({u}^n_j)_{j \in \mathbb{Z}}$ at time $t^n$, it is thus a question of
defining its evolution towards the next time level $t^{n+1}$. More
precisely, and in the context of a finite volume conservative scheme,
we have to define the numerical fluxes $(f^n_{j+1/2})_{j \in
\mathbb{Z}}$ coming in (\ref{schema_cons}). For that, we still assume 
\begin{equation} \label{paspointsonique}
\,\, \mbox{either}\,\, 
f'(u) \geq 0 \,\, \, \mbox{for all } u,
 \quad
\mbox{or}\,\, f'(u) \leq 0 \,\, \, \mbox{for all } u,
\end{equation}
so that propagation is in one direction only. 
According to the previous section, information in cell 
$\mathcal{C}_j$ is 
understood as an element of the inner structure of the Riemann problem
associated with initial states $u^n_{j-1}$ and $u^n_{j+1}$. This one
will be used to compute either $f^n_{j+1/2}$ (if $f'(u) \geq 0$) or
$f^n_{j-1/2}$ (if $f'(u) \leq 0$). \\ 
In Section~\ref{2-0}, it is stated that the Riemann problem 
associated with initial states $u^n_{j-1}$ and $u^n_{j+1}$ may contain
a nonclassical shock between $u^n_{j-1}$ and
$\varphi^{\flat}(u^n_{j-1})$ if the function is concave-convex (and
between $\varphi^{-\flat}(u^n_{j+1})$ and $u^n_{j+1}$ if the function
is convex-concave). 

Recall that these nonclassical waves are difficult to capture numerically 
and require special attention.  
(We have shown in the previous section that as many others, the upwind scheme does not suit.) 
Instead of considering $u^n_j$ as a sufficiently accurate information for the
structure of the Riemann solution associated with the initial states
$u^n_{j-1}$ and $u^n_{j+1}$, we propose to replace it (whenever possible) with a discontinuity 
separating 
$u^n_{j,l} = \varphi^{-\flat}(u^n_{j+1})$ on the left and $u^n_{j,r} =
\varphi^{\flat}(u^n_{j-1})$ on the right, and located at point
$\overline{x}_j \in \mathcal{C}_j$. In other words, we propose to
introduce in the cell $\mathcal{C}_j$ the right (respectively left)
state $\varphi^{\flat}(u^n_{j-1})$ (respectively
$\varphi^{-\flat}(u^n_{j+1})$) of the nonclassical discontinuity which
is expected to be present in the Riemann solution associated with
$u^n_{j-1}$ and $u^n_{j+1}$ (depending on if $f$ obeys
(\ref{conca-conve}) or (\ref{conve-conca})). As in the previous
section, one requires the reconstructed discontinuity to satisfy the
conservation property (\ref{x_disc}) and to be located inside
$\mathcal{C}_j$, that is $0 \leq d^n_j \leq 1$ with $d^n_j$ given in
(\ref{defdnj}).  
Here, we
let $u^n_{j,l} = u^n_{j,r} = u^n_j$ if $d^n_j$ given in (\ref{defdnj})
does not belong to $[0,1]$.

\noindent Then, we naturally set for all $j \in \mathbb{Z}$: \\
$(i)$ if $f$ is non-decreasing 
\begin{equation} \label{deffluxfpp}
\Delta t f^n_{j+1/2} = 
\left\{
\begin{array}{ccl}
\min(\Delta t_{j+1/2}, \Delta t) f(u^n_{j,r}) + & 
\max(\Delta t - \Delta t_{j+1/2}, 0) f(u^n_{j,l}), 
\\ 
& 0 \leq d^n_j \leq 1,  \\
\Delta t f(u^n_j), & 
\\
&\mbox{otherwise}, 
\end{array}
\right.
\end{equation}
with 
\begin{equation} \label{deltat_int_3a}
\Delta t_{j+1/2} = \frac{1 - d^n_j}{\sigma(u^n_{j,l},u^n_{j,r})}
\Delta x. 
\end{equation}
$(ii)$ if $f$ is non-increasing:
\begin{equation} \label{deffluxfpn}
\Delta t f^n_{j-1/2} = 
\left\{
\begin{array}{ccl}
\min(\Delta t_{j-1/2}, \Delta t) f(u^n_{j,l}) + & 
\max(\Delta t - \Delta t_{j-1/2}, 0) f(u^n_{j,r}), 
\\
&  
0 \leq d^n_j \leq 1, \\
\Delta t f(u^n_{j}), &
\\
& \mbox{otherwise}, 
\end{array}
\right.
\end{equation}
with 
\begin{equation} \label{deltat_int_3b}
\Delta t_{j-1/2} = \frac{d^n_j}{-\sigma(u^n_{j,l},u^n_{j,r})} \Delta
x. 
\end{equation}
Note that contrary to the linear advection (see the first
illustration in the previous section), the local time step $\Delta
t_{j+1/2}$ (respectively $\Delta t_{j-1/2}$) given by
(\ref{deltat_int_3a}) (respectively (\ref{deltat_int_3b})) is now only
a prediction of the time needed by the reconstructed discontinuity to
reach the interface $x_{j+1/2}$ (respectively $x_{j-1/2}$).  
The prediction step is however exact in the case of an isolated
nonclassical discontinuity (see the second illustration in the
previous section) and more generally as soon as $u^n_{j-1}$ and
$u^n_{j+1}$ verify $u^n_{j+1} = \varphi^{\flat}(u^n_{j-1})$. 
\ \\

Observe that the proposed scheme belongs to the class of five-point
schemes, since ${u}^{n+1}_j$ depends on ${u}^{n}_{j-2}$,
${u}^{n}_{j-1}$, ${u}^{n}_j$, ${u}^{n}_{j+1}$ and ${u}^{n}_{j+2}$. 


\subsection*{Stability and consistency properties}

\noindent We now state and prove important properties enjoyed by our algorithm.

We assume that the flux $f$ satisfies the
monotonicity condition (\ref{paspointsonique}) and either the
concave-convex or concave-convex conditions (\ref{conca-conve}) or
(\ref{conve-conca}) respectively. Then, under the CFL restriction 
\begin{equation} 
\label{cfl}
\frac{\Delta t}{\Delta x} \max_{u} |f'(u)| \leq 1,
\end{equation}
where the maximum is taken over all the $u$ under consideration, the
conservative scheme (\ref{schema_cons}) with $f^n_{j+1/2}$ defined 
for all $j \in \mathbb{Z}$ by (\ref{deffluxfpp})-(\ref{deffluxfpn}) is
consistent with (\ref{1})-(\ref{2})-(\ref{kc}) in the following sense. 

\begin{property}[Flux consistency.] 
\label{prop1} 
Assume that
${u}:={u}^{n}_{j-1} = {u}^{n}_j = {u}^{n}_{j+1}$, then
${f}^{n}_{j+1/2} = f({u})$ if $f'\geq 0$ and 
${f}^{n}_{j-1/2} = f({u})$ if $f' \leq 0$. 
\end{property}

\begin{property}[Classical solutions.] 
\label{prop2}
Assume that
${u}^{n}_{j-2}$, ${u}^{n}_{j-1}$, ${u}^{n}_j$, ${u}^{n}_{j+1}$ and
${u}^{n}_{j+2}$ belong to the same region of convexity of $f$. Then
the definition ${u}^{n+1}_{j}$ given by the conservative scheme
(\ref{schema_cons})-(\ref{deffluxfpp})-(\ref{deffluxfpn}) coincides
with the one given by the usual upwind conservative scheme. Then it
obeys all the usual stability properties provided by this scheme. In
particular, the strategy is convergent if the whole discrete solution
belongs to the same region of convexity of $f$. 
\end{property}

\begin{property}[Isolated nonclassical shock waves.] 
\label{prop3}
Let $u_l$ and $u_r$ be two initial states such that 
$u_r = \varphi^{\flat}(u_l)$. Assume that ${u}^{0}_{j} = {u}_l$ if $j
\leq 0$ and ${u}^{0}_{j} = {u}_r$ if $j \geq 1$. Then the conservative
scheme (\ref{schema_cons})-(\ref{deffluxfpp})-(\ref{deffluxfpn})
provides an exact numerical solution on each cell $\mathcal{C}_j$ in
the sense that 
\begin{equation} \label{seit}
u^n_j = \frac{1}{\Delta x} \int_{x_{j-1/2}}^{x_{j+1/2}} u(x,t^n) dx,
\,\, \qquad \, j \, \in \, \mathbb{Z}, \,\, \, n \, \in \, \mathbb{N}, 
\end{equation}
where $u$ denotes the exact Riemann solution
of (\ref{1})-(\ref{2})-(\ref{ditr})-(\ref{kc}) given by
$u(x,t) = u_l$ if $x < \sigma(u_l,u_r) t$ and $u(x,t) = u_r$
otherwise, and is convergent towards $u$. In particular, the numerical
discontinuity is diffused on one cell at most.
\end{property}

The following comments are in order. Property $(i)$ shows
that the proposed numerical flux function is consistent in the
classical sense of finite volume methods. 
Properties $(ii)$ and $(iii)$ provide us with crucial stability/accuracy properties. 
They show that the method is actually convergent
if the solution remains in the same
region of convexity of $f$ (see $(ii)$) or, more importantly, the solution consists in
an isolated nonclassical discontinuity satisfying the prescribed
kinetic relation (see $(iii)$). 
We emphasize that all of the {\sl conservative} schemes proposed so far in the literature 
violate the latter property.

\vskip.2cm 

\noindent{\bf Proof of Property~\ref{prop1}.}  
(i) If ${u}:={u}^{n}_{j-1} = {u}^{n}_j = {u}^{n}_{j+1}$ then 
$$
d^n_j = \frac{\varphi^{\flat}(u) - u}{\varphi^{\flat}(u) -
\varphi^{-\flat}(u)}. 
$$
The property $0 \leq d^n_j \leq 1$ means 
$\min(\varphi^{-\flat}(u),\varphi^{\flat}(u)) \leq u \leq 
\max(\varphi^{-\flat}(u),\varphi^{\flat}(u))$ and cannot hold, 
since $u$ and $\varphi^{\flat}(u)$ do not have the same sign for all
$u$. Then, we obtain ${f}^{n}_{j+1/2} = f({u})$ if $f'\geq 0$ and
${f}^{n}_{j-1/2} = f({u})$ if $f' \leq 0$ by
(\ref{deffluxfpp})-(\ref{deffluxfpn}).

\

\noindent{\bf Proof of Property~\ref{prop2}.} 
Assume without restriction that 
$f' \geq 0$ and recall that $0 \leq d^n_{j-1} \leq 1$ 
and $0 \leq d^n_{j} \leq 1$ respectively means that
$$
\min(\varphi^{-\flat}({u}^{n}_{j}),\varphi^{\flat}({u}^{n}_{j-2}))
\leq {u}^{n}_{j-1} \leq
\max(\varphi^{-\flat}({u}^{n}_{j}),\varphi^{\flat}({u}^{n}_{j-2}))
$$
and
$$
\min(\varphi^{-\flat}({u}^{n}_{j+1}),\varphi^{\flat}({u}^{n}_{j-1}))
\leq {u}^{n}_{j} \leq
\max(\varphi^{-\flat}({u}^{n}_{j+1}),\varphi^{\flat}({u}^{n}_{j-1})).
$$
These inequalities are not valid since by definition $u$ and
$\varphi^{\flat}(u)$ do not belong to the same region of convexity of
$f$. By (\ref{deffluxfpp})-(\ref{deffluxfpn}), the numerical fluxes
$f^n_{j \pm 1/2}$ coincides with the usual upwind fluxes and the
conclusion follows.

\

\noindent{\bf Proof of Property~\ref{prop3}.} 
First, note that there is no relevant reconstruction in the
first iteration. Indeed, the property $0 \leq d^n_j \leq 1$ reads as
follows if $j<0$ or $j>1$, 
$$
\begin{array}{ccc}
& 0 \leq d^n_j \leq 1 \mbox{ if and only if } 
\\
& \left\{
\begin{array}{ccc}
\min(\varphi^{-\flat}(u_l),\varphi^{\flat}(u_l)) \leq u_l \leq 
\max(\varphi^{-\flat}(u_l),\varphi^{\flat}(u_l)), &  & j < 0,
\\ 
\min(\varphi^{-\flat}(u_r),\varphi^{\flat}(u_r)) \leq u_r \leq 
\max(\varphi^{-\flat}(u_r),\varphi^{\flat}(u_r)), &  & j > 1,
\end{array}
\right. 
\end{array}
$$
which again cannot hold (see $(i)$ below), while if $j=0$ or $j=1$,
the relation $u_r=\varphi^{\flat}(u_l)$ and definition (\ref{defdnj})
give 
$$
\left\{
\begin{array}{rcl}
d^n_j = {\displaystyle{\frac{u_r-u_l}{u_r-u_l}}} = 1, &  & j =
0, \\ 
d^n_j = {\displaystyle{\frac{u_r-u_r}{u_r-u_l}}} = 0, &  & j =
1, \\ 
\end{array}
\right.
$$
so that the reconstructions exist but are trivial:
$u_l=\varphi^{-\flat}(u_r)$ (respectively $u_r=\varphi^{\flat}(u_l)$)
takes the whole cell associated with $j=0$ (respectively $j=1$). \\ 
Assume now without restriction that $f$ is non-decreasing and let
$\Delta t$ be such that (\ref{cfl}) holds. After one time step $\Delta
t$, the exact solution given by $u(x,\Delta t) = u_l$ if $x <
\sigma(u_l,u_r) \Delta t$ and $u(x,\Delta t) = u_r$ otherwise is such
that 
\begin{equation} \label{seitpr}
\frac{1}{\Delta x} \int_{x_{j-1/2}}^{x_{j+1/2}} u(x,\Delta t) dx
=
\left\{
\begin{array}{ccc}
u_l, &  & j \leq 0, \\
u_r - \sigma(u_l,u_r) \frac{\Delta t}{\Delta x} (u_r - u_l),& 
& j = 1, \\ 
u_r, &  & j > 1. \\
\end{array}
\right.
\end{equation}
But recall that $\sigma(u_l,u_r)$ is given by (\ref{srh}) so that we
have 
\begin{equation} \label{seitpr2}
\frac{1}{\Delta x} \int_{x_{j-1/2}}^{x_{j+1/2}} u(x,\Delta t) dx
=
\left\{
\begin{array}{ccc}
u_l - \frac{\Delta t}{\Delta x} (f(u_l) - f(u_l)), &  & j \leq
0, \\ 
u_r - \frac{\Delta t}{\Delta x} (f(u_r) - f(u_l)),& & j = 1,
\\ 
u_r - \frac{\Delta t}{\Delta x} (f(u_r) - f(u_r)), &  & j > 1,
\\ 
\end{array}
\right.
\end{equation}
that is
\begin{equation} \label{cppit}
 u^1_j = \frac{1}{\Delta x} \int_{x_{j-1/2}}^{x_{j+1/2}} u(x,\Delta t)
dx, \,\, \qquad \, j \, \in \, \mathbb{Z}.
\end{equation}
The identity (\ref{seit}) is then proved for the first iterate. \\
What happens now in the next time iteration ? At this stage, it is
first clear (see the previous discussion just below) that only cell
$\mathcal{C}_1$ is going to be dealt with a reconstruction. Now, the
main point of the proof lies in the fact that the reconstructed
discontinuity in this cell actually joins the expected states
$\varphi^{-\flat}(u^1_2) = \varphi^{-\flat}(u_r) = u_l$ and
$\varphi^{\flat}(u^1_0) = \varphi^{\flat}(u_l) = u_r$ and is located
exactly at point $x = \sigma(u_l,u_r) \Delta t$ by the conservation
property (\ref{cppit}). In other words, we have reconstructed the
exact solution at time $t=\Delta t$. To derive the required identity
(\ref{seit}) for the second iterate, it is sufficient to recall that
by Green's formula the conservative scheme (\ref{schema_cons}) with
$f^n_{j+1/2}$ defined for all $j \in \mathbb{Z}$ by
(\ref{deffluxfpp})-(\ref{deffluxfpn}) is equivalent for $n=2$ to
average the evolution of this exact solution up to time $t^2 = 2
\Delta t$. And the process is going on in a similar way for the next
time iterations, which proves the result.


\section{Numerical experiments}

We mostly consider here the flux $f(u) := u^3 + u$, thus $f$ is concave-convex in the sense given in the
second section. For the entropy-entropy flux pair $(U,F)$ required in (\ref{2}), we use 
$$
U(u) := u^2, \quad F(u) := \frac{3}{2} u^4 + u^2.
$$
Easy calculations lead to explicit formulas for $\varphi^{\natural}$
and $\varphi^{-\natural}$:
$$
\varphi^{\natural} = -\frac{u}{2}, 
\qquad 
\varphi^{-\natural} = -2u, 
\qquad 
\varphi^{\flat}_0(u) = -u. 
$$
Moreover, we have here $\varphi^{\sharp}(u) = - u - \varphi^{\flat}(u)$. \\

The choice of the kinetic function $\varphi^{\flat}$ must be in
agreement with relations (\ref{hypkf1}) with $\varphi^{\natural}$ and
$\varphi^{\flat}_0$ just calculated. Here, we will choose the kinetic function 
$$
\varphi^{\flat}(u) = -\beta u, \quad \beta \in \left[0.5,1\right),
$$
which, as observed in Bedjaoui and LeFloch \cite{bl}, 
can be realized by an augmented model based on nonlinear diffusion and dispersion terms. 
In the following, we will take $\beta=0.75$. \\ 
\ \\

\begin{figure}[htbp]
\begin{center}
\includegraphics[width=8.0cm]{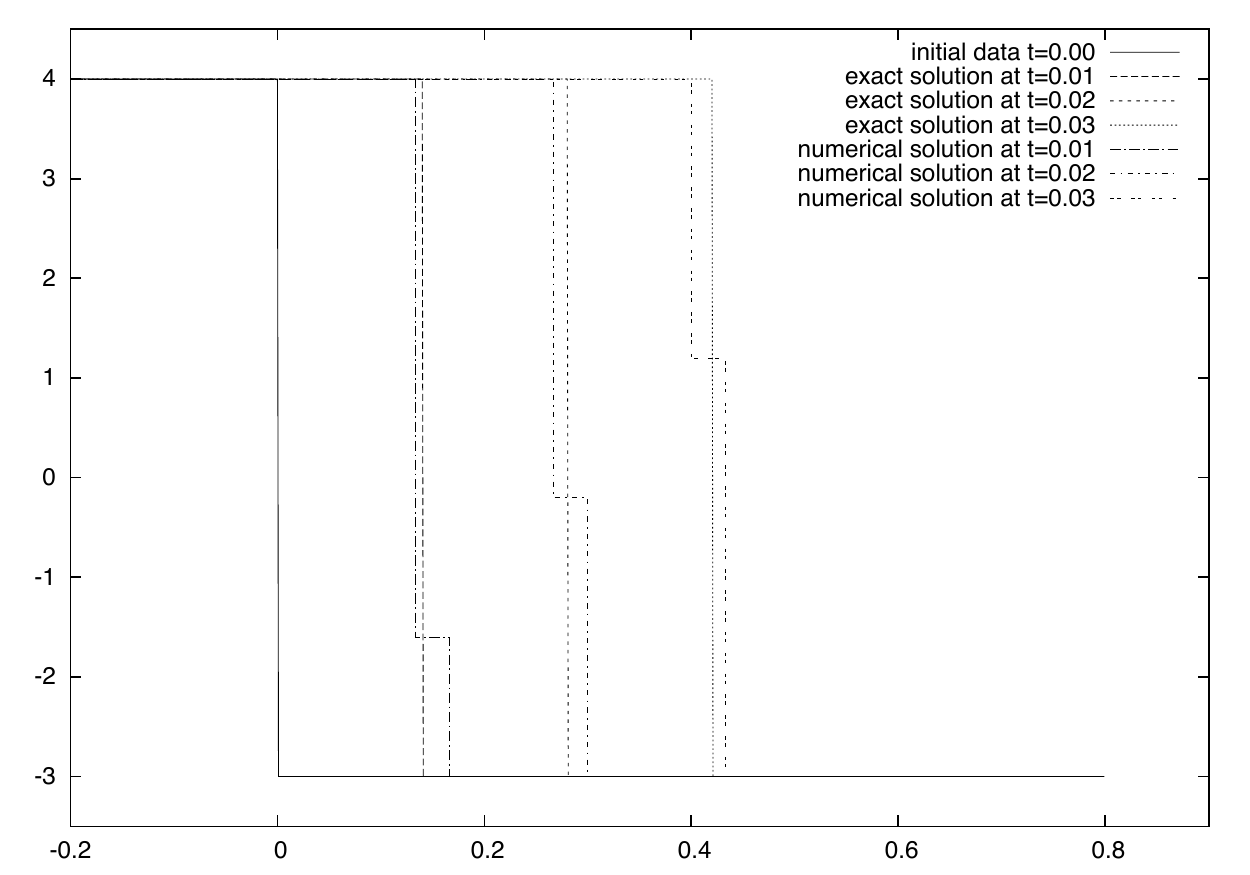}
\vspace{-0.45cm}
\caption{Test A - Nonclassical shock -- 30 points} 
\label{Test0}
\end{center}
\end{figure}

{\bf Test A.} Let us check Property~\ref{prop3} numerically, 
which is concerned with the exact capture of isolated nonclassical shocks. Thus, consider the following
nonclassical shock as a Riemann initial condition 
$$
u_0(x)=
\left\{
\begin{array}{rl}
4, & x<0,\\
\varphi^\flat(4)=-3, & x>0,
\end{array}
\right.
$$
The numerical solution shown in Figure \ref{Test0} is exact everywhere but 
in the single cell containing the nonclassical shock. (We
sometimes use a piecewise constant representation in the figure, in order
for the interpretation of the numerical solutions to be easier.) However, 
as expected, the value in this cell 
coincides with the average of the corresponding exact solution (see (\ref{seit})),
and allows (after reconstruction) to recover the {\sl exact} location of the
discontinuity (using the conservation property of scheme).
This property explains why the numerical solution stays sharp when the time evolves. 
\\

\begin{figure}[htbp]
\begin{center}
\includegraphics[width=8.0cm]{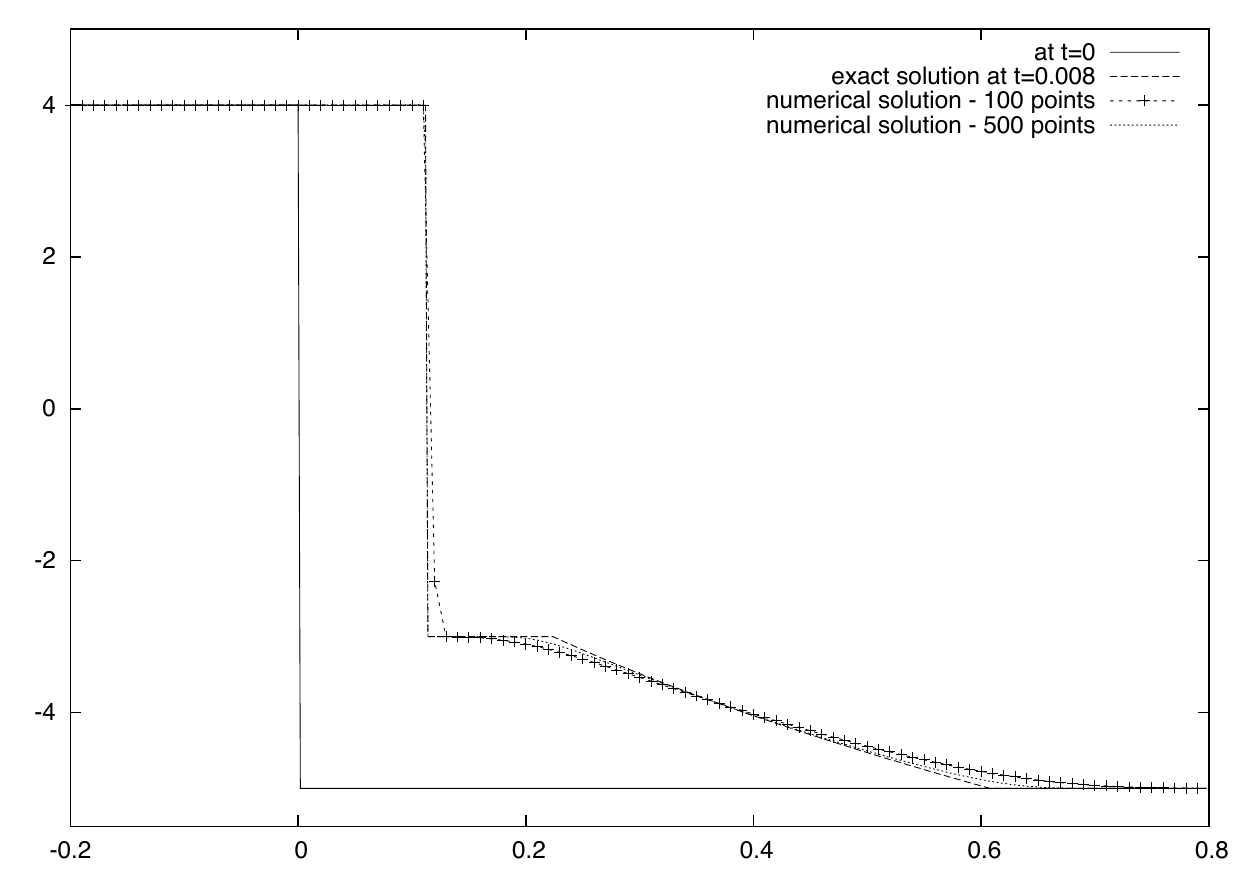}
\includegraphics[width=8.0cm]{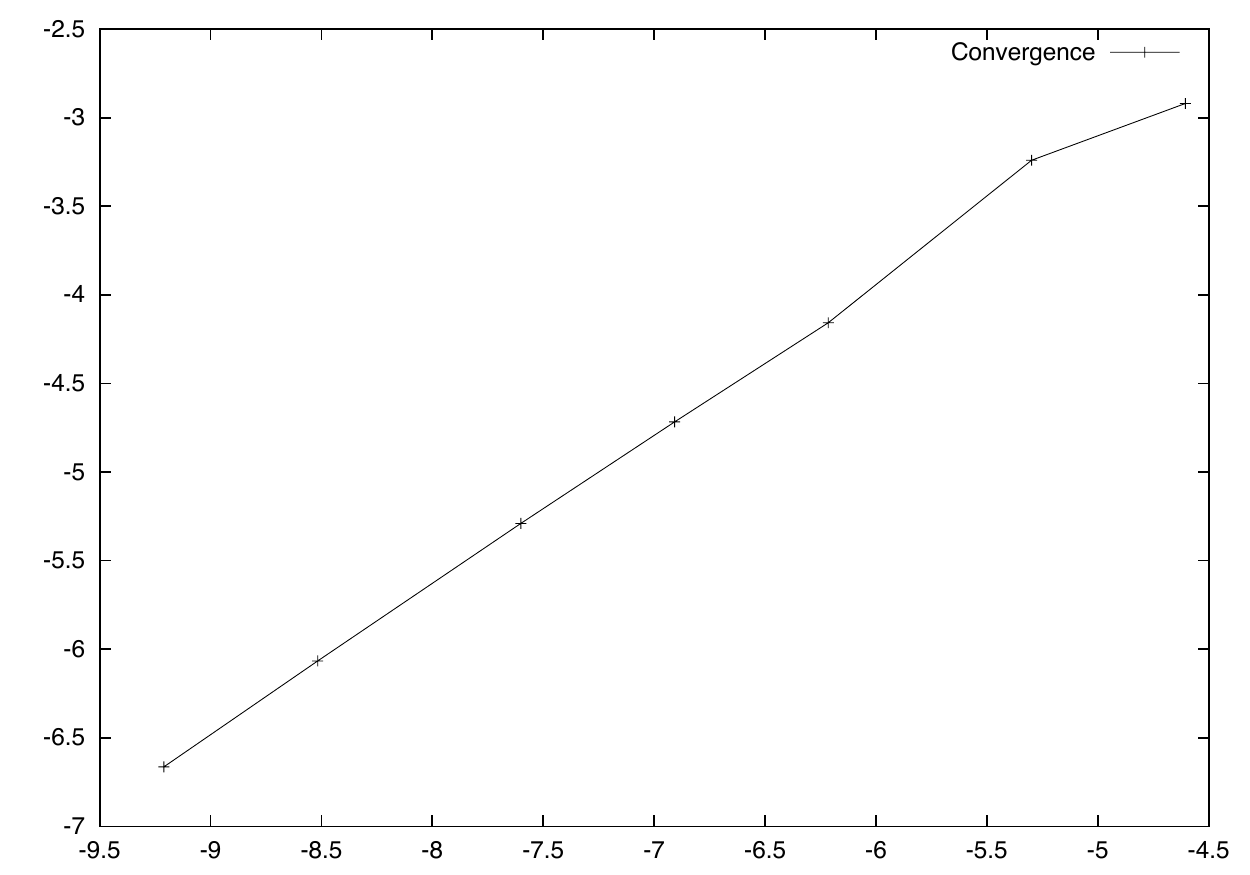}
\vspace{-0.45cm}
\caption{Test B - Nonclassical shock and rarefaction -- $L^1$ convergence
 (log($E_{L^1}$) versus log($\Delta x$))} 
\label{Test1}
\end{center}
\end{figure}

{\bf Test B.} In our second test we consider the Riemann problem with initial data
$$
u_0(x)=
\left\{
\begin{array}{rl}
4, & x<0,\\
-5, & x>0,
\end{array}
\right.
$$
whose solution is a nonclassical shock followed by a rarefaction
wave. The two left-hand curves in Figure \ref{Test1} are
performed with $\Delta x=0.01$ and $\Delta x=0.002$, respectively.
The nonclassical shock, as previously, is localized in
a single computational cell.\\
The right-hand figure represents the logarithm of the $L^1$-error 
(between the exact and the numerical solution) versus the logarithm of $\Delta x$. 
The numerical order of convergence is about 0.8374.\\

\begin{figure}[htbp]
\begin{center}
\includegraphics[width=8.0cm]{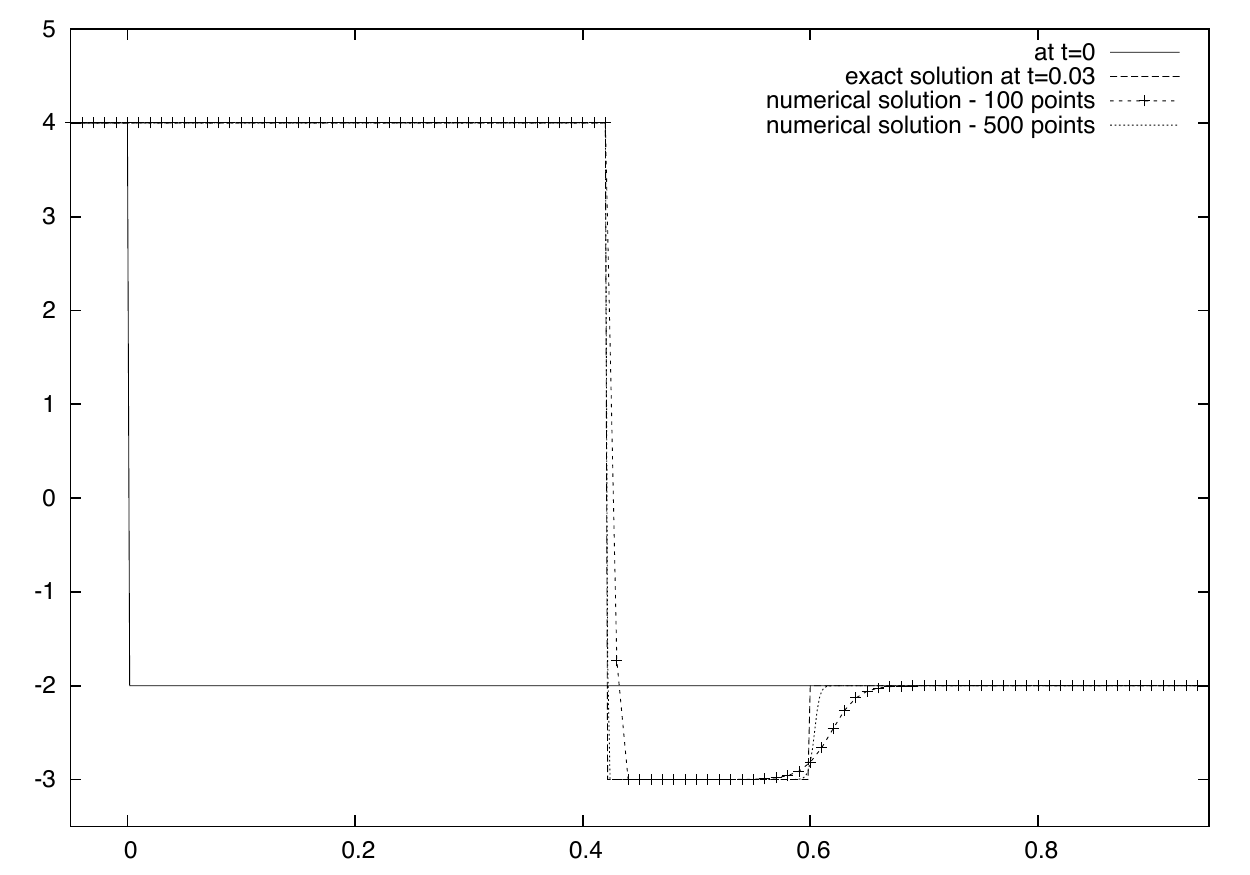}
\includegraphics[width=8.0cm]{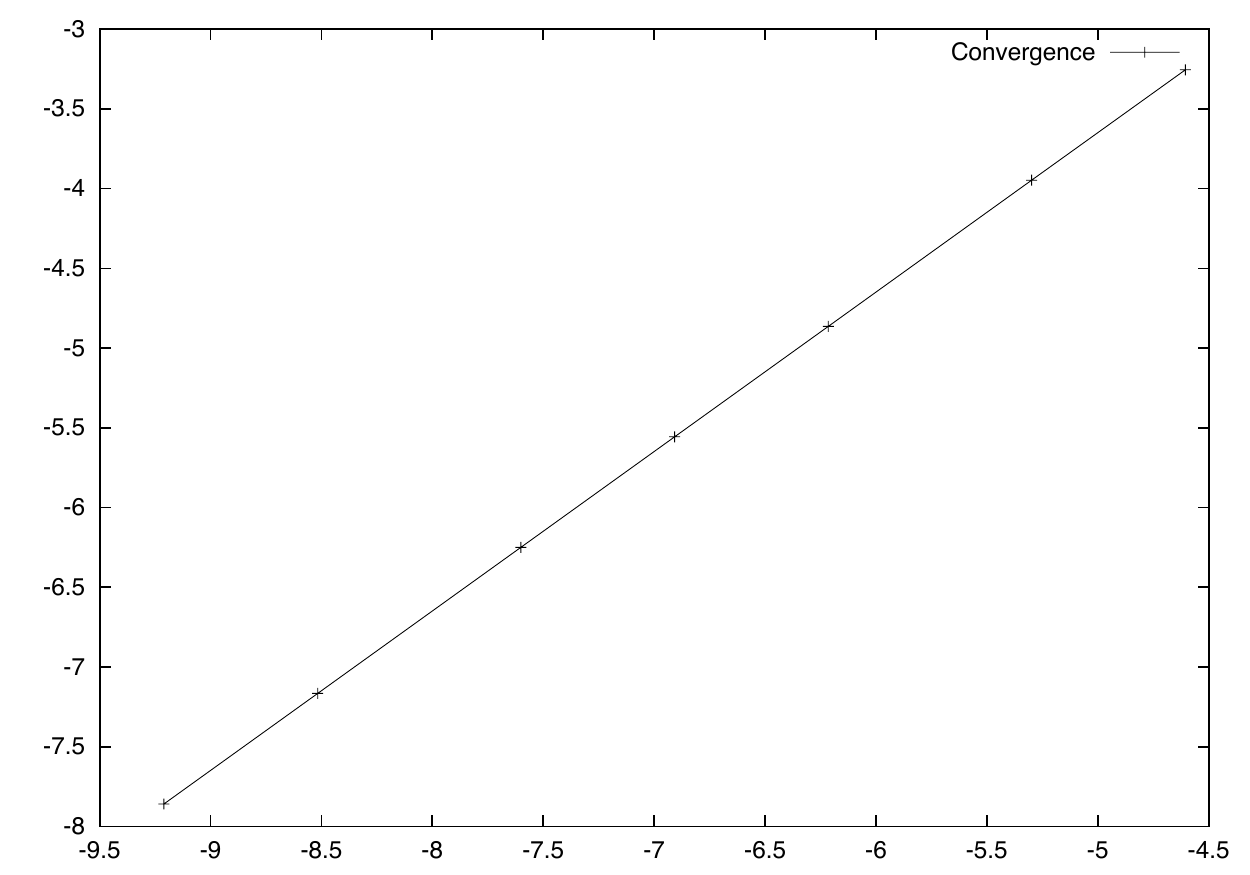}
\vspace{-0.45cm}
\caption{Test C - Nonclassical and classical shocks -- $L^1$ convergence
(log($E_{L^1}$) versus log($\Delta x$))} 
\label{Test2}
\end{center}
\end{figure}

{\bf Test C} (Figure \ref{Test2}). Now, we choose another Riemann initial
condition which develops a nonclassical shock followed by a classical shock:  
$$
u_0(x)=
\left\{
\begin{array}{rl}
4, & x<0,\\
-2, & x>0.
\end{array}
\right.
$$
We can make the same observation as previously, concerning the
nonclassical shock; it is sharply captured and arises in a small spatial domain.
However, note here that the classical shock does contain some numerical diffusion: in fact, 
our scheme is exactly the upwinding scheme if the values of the solution 
remains in a given convexity region for the flux $f$.\\
Once again, the plot with the $L^1$-error shows the numerical convergence with order about 0.9999.\\

\begin{figure}[htbp]
\begin{center}
\includegraphics[width=8.0cm]{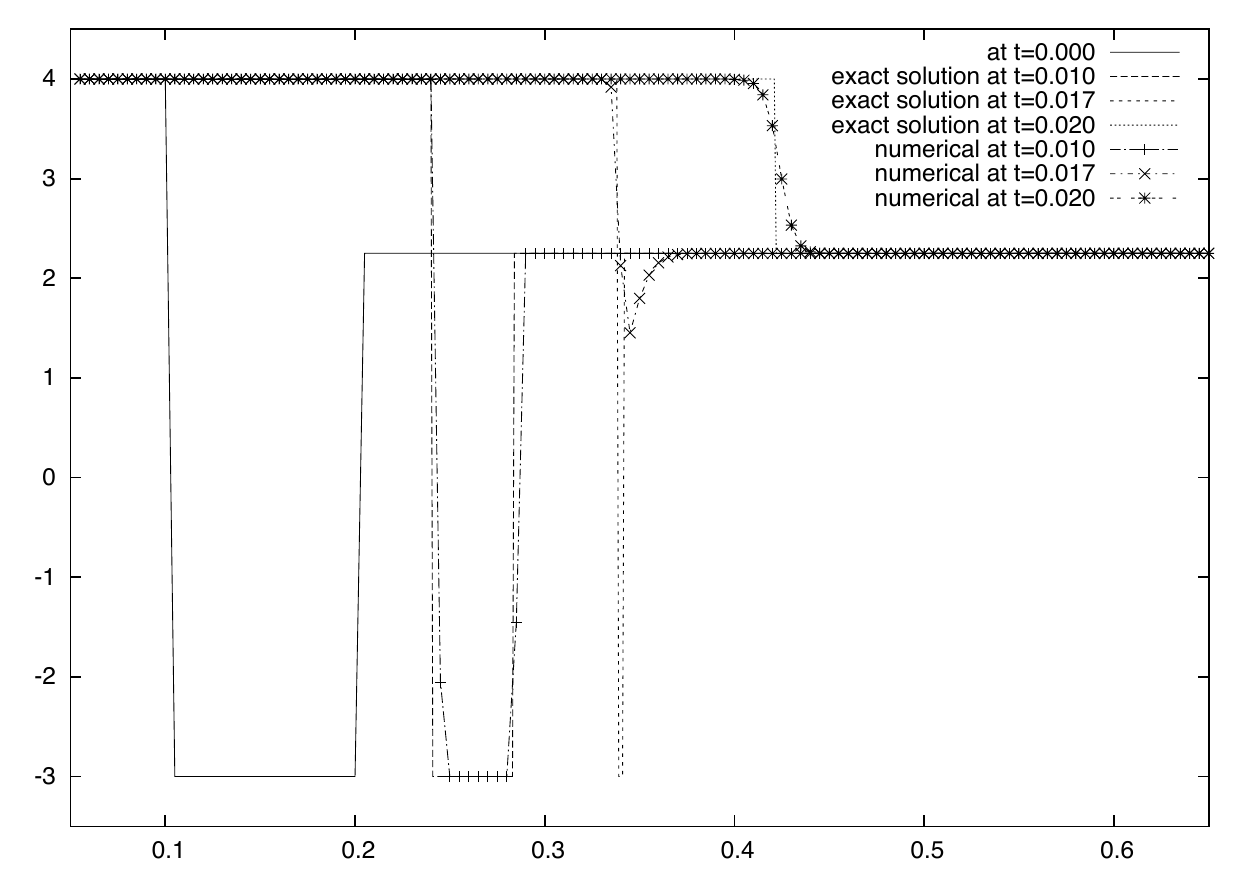}
\vspace{-0.45cm}
\caption{Test D - Interaction of two nonclassical shocks}
\label{Test3}
\end{center}
\end{figure}

{\bf Test D} (Figure \ref{Test3}). We now take an initial data composed of two nonclassical shocks 
that interact: 
$$
u_0(x)=
\left\{
\begin{array}{ll}
4 = \varphi^{-\flat}(-3), & x<0.1\\
-3, & 0.1<x<0.2\\
2.25= \varphi^{\flat}(-3), & x>0.2.
\end{array}
\right.
$$
The computation is performed with $\Delta x=0.05$ and plotted at four
successive times $t=0, 0.0010, 0.0017$, and $0.0020$. We observe that
the two nonclassical shocks cancel each other at the interaction, and generate a single
classical shock, in accordance with the general theory in \cite{lefloch1}.

\begin{figure}[htbp]
\begin{center}
\includegraphics[width=8.0cm]{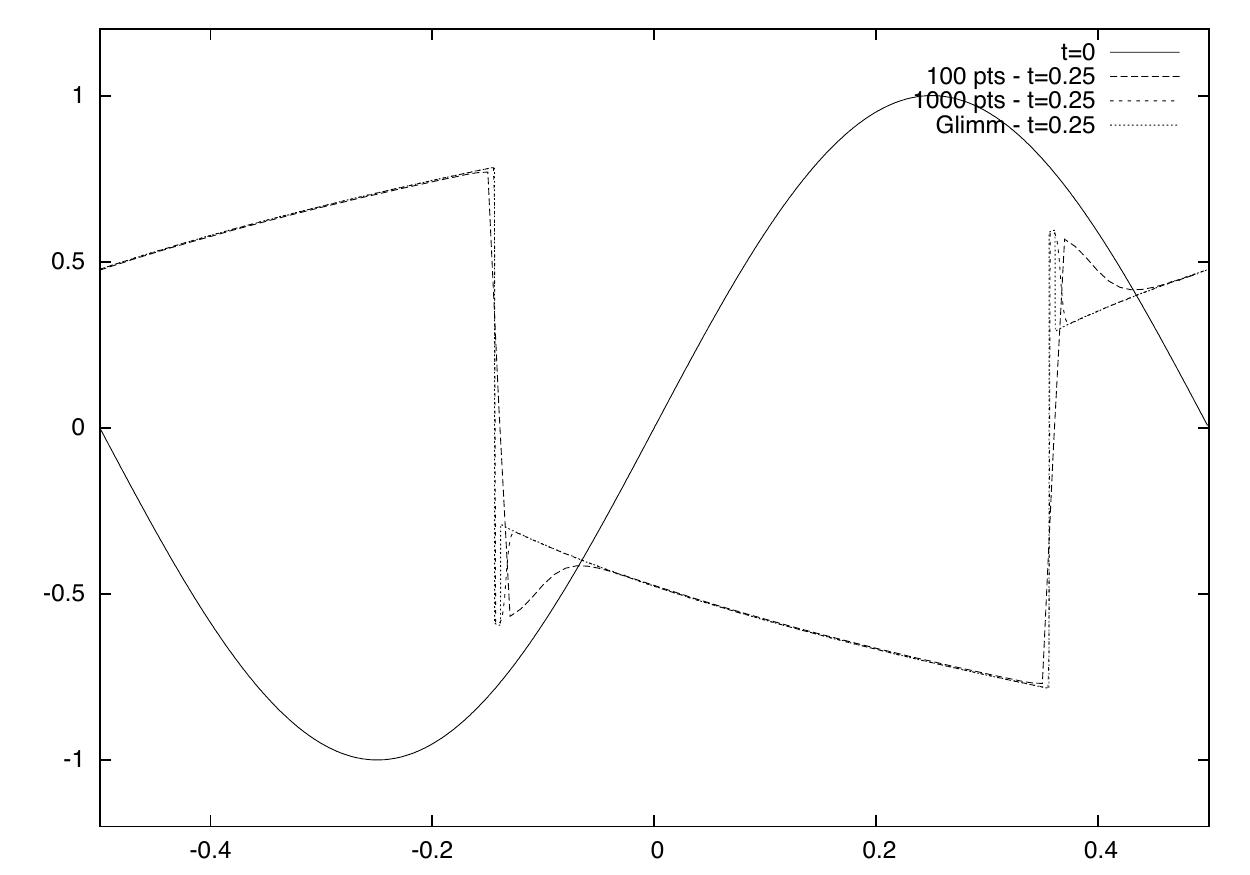}
\includegraphics[width=8.0cm]{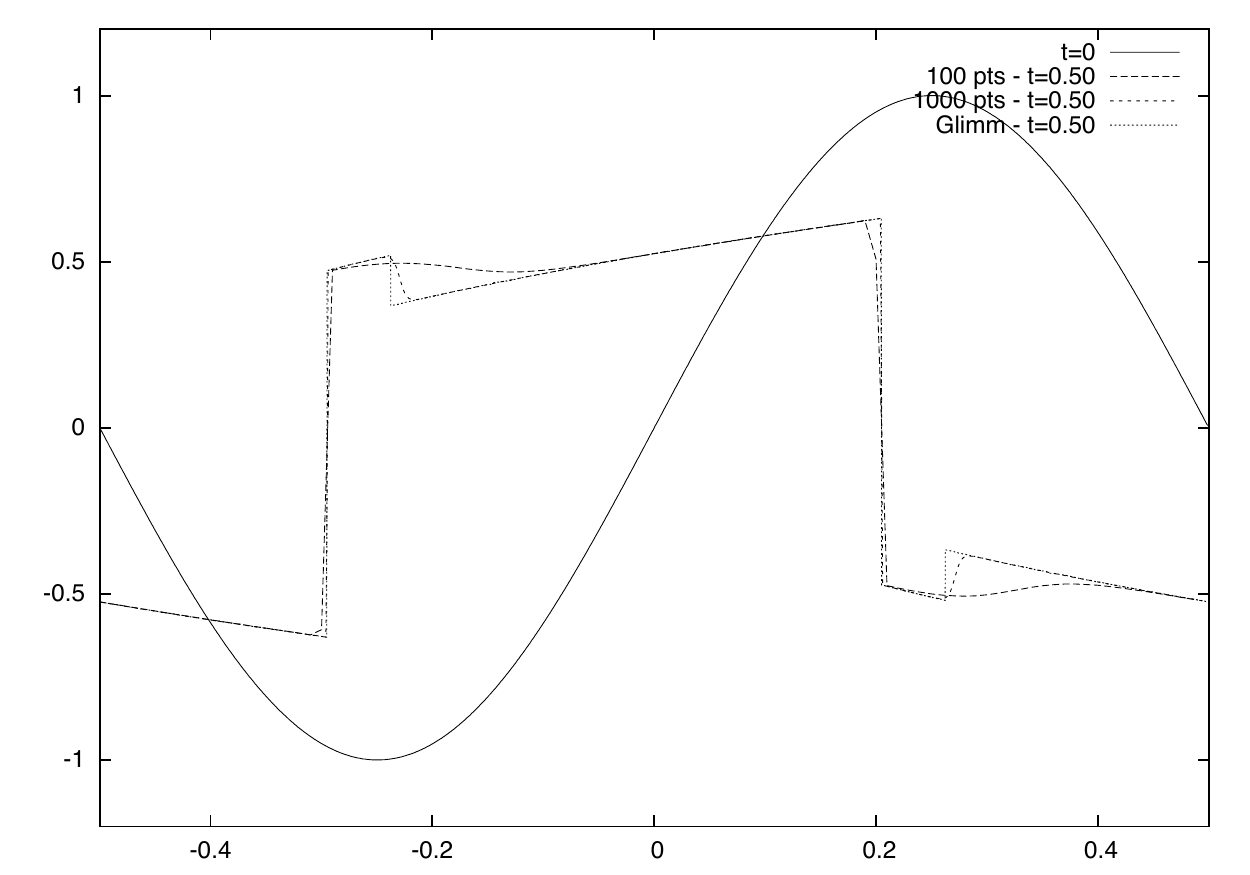}
\vspace{-0.45cm}
\caption{Test E - Periodic initial data - reconstruction scheme and
 Glimm scheme} 
\label{Test4}
\end{center}
\end{figure}

{\bf Test E} (Figure \ref{Test4}). Next, we consider the periodic initial
condition 
$$
u_0(x)=\sin\left({x\over 2\pi}\right),
$$
with periodic boundary conditions $u(-0.5,t)=u(0.5,t)$. The exact
solution is not known explicitly, so we compare our numerical solution
with the solution generated by Glimm's random choice scheme \cite{glimm} in which  
we have replaced the classical solver by the nonclassical solver described in Section~\ref{2-0}. 
We use here van der Corput's random sequence $(a_n)$, defined by 
$$
a_n=\sum_{k=0}^m i_k 2^{-(k+1)},
$$
where $n=\sum_{k=0}^m i_k 2^k,\ i_k\in\{0,1\}$, denotes the binary
expansion of the integer $n$. Figure \ref{Test4} represents the
solutions at the times $t=0, 0.25$ and $0.5$ for our scheme with $\Delta
x=0.01$ and with $\Delta x=0.0001$,
and for the Glimm scheme with $\Delta x=0.0001$ (to serve as a reference). 
The two methods strongly agree. Roughly speaking, the increasing parts of $u_0$ 
evolve as rarefactions,
while the decreasing parts are compressed and develop in a classical shock and, then, 
when left- and right-hand states at the shocks change sign, nonclassical shocks (which 
do satisfy the expected kinetic relation) and new faster classical shocks on the right-hand side arise. \\

\begin{figure}[htbp]
\begin{center}
\includegraphics[width=8.0cm]{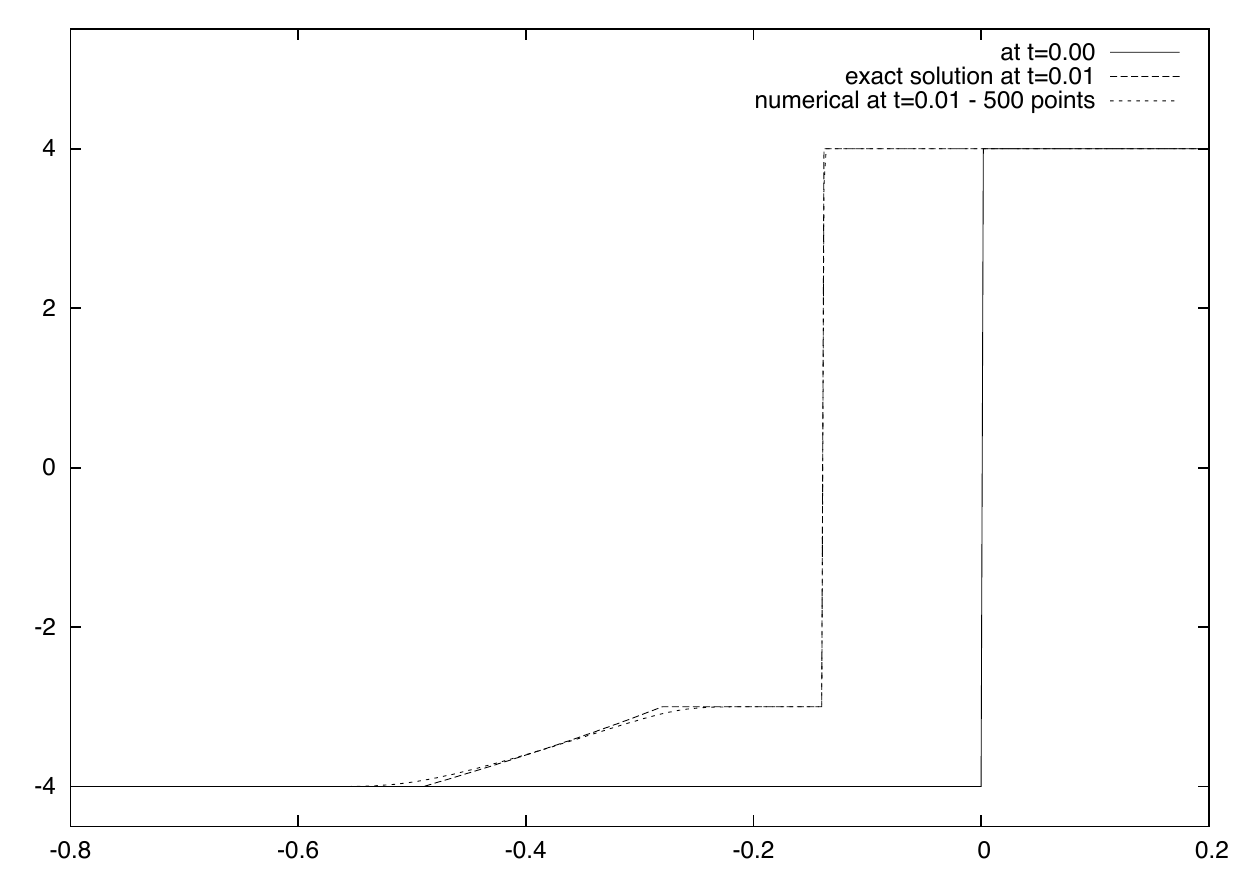}
\includegraphics[width=8.0cm]{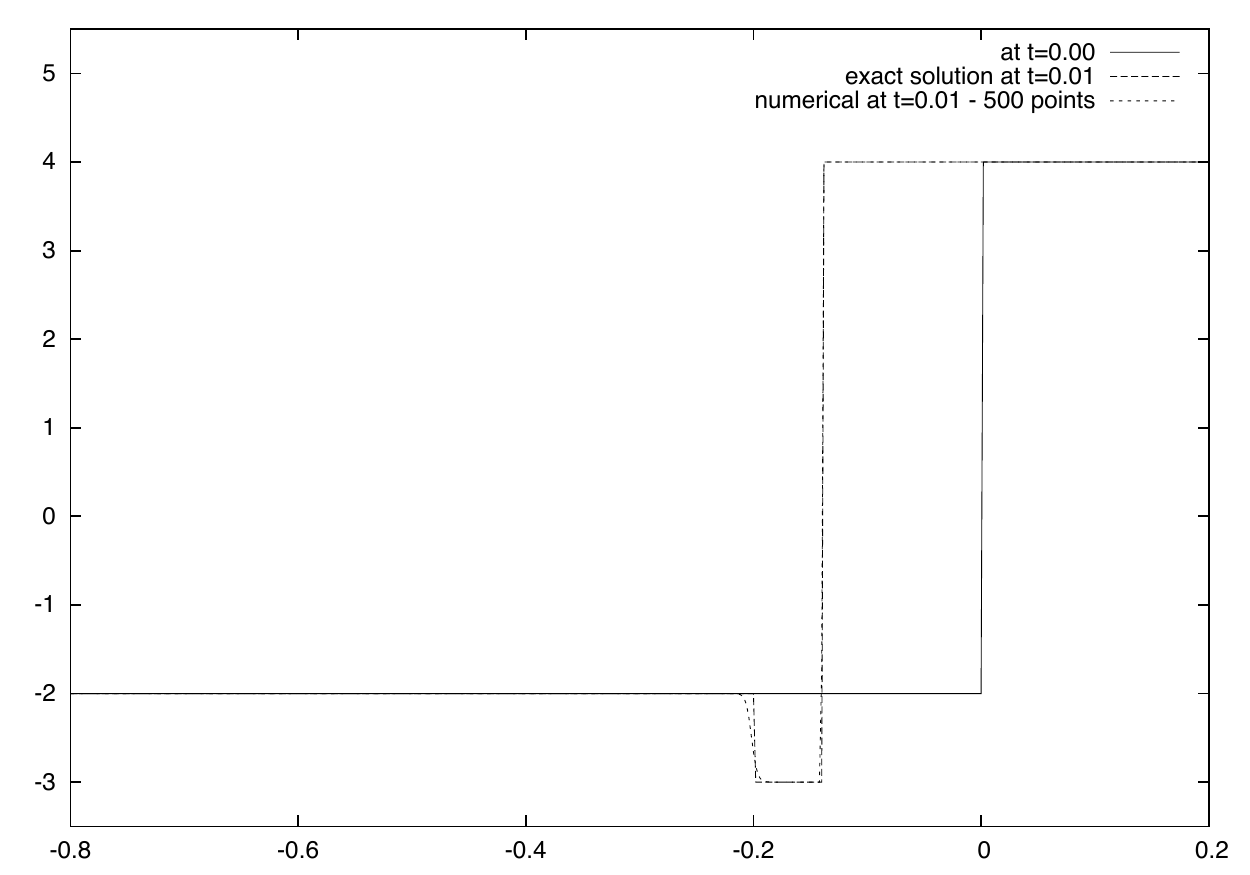}
\vspace{-0.45cm}
\caption{Test E - Two examples in the convex-concave case} 
\label{Test5}
\end{center}
\end{figure}

{\bf Test F} (Figure \ref{Test5}). To illustrate the behavior of convex-concave flux functions, 
we finally compute two Riemann solutions with opposite flux $f(u)=-u^3-u$ (so $f'<0$ and the solutions move from
right to left) and the same kinetic function $\varphi^\flat(u)=-0.75\
u$: the first one (left-hand figure) corresponds to the initial data
$$
u_0(x)=
\left\{
\begin{array}{rl}
-4, & x<0,\\
4, & x>0,
\end{array}
\right.
$$
and develops a rarefaction wave and a nonclassical shock; the second
one (right-hand figure) corresponds to the initial data
$$
u_0(x)=
\left\{
\begin{array}{rl}
-2, & x<0,\\
4, & x>0,
\end{array}
\right.
$$
and the corresponding solution is a classical shock followed by a
nonclassical shock.\\

\begin{figure}[htbp]
\begin{center}
\includegraphics[width=8.0cm]{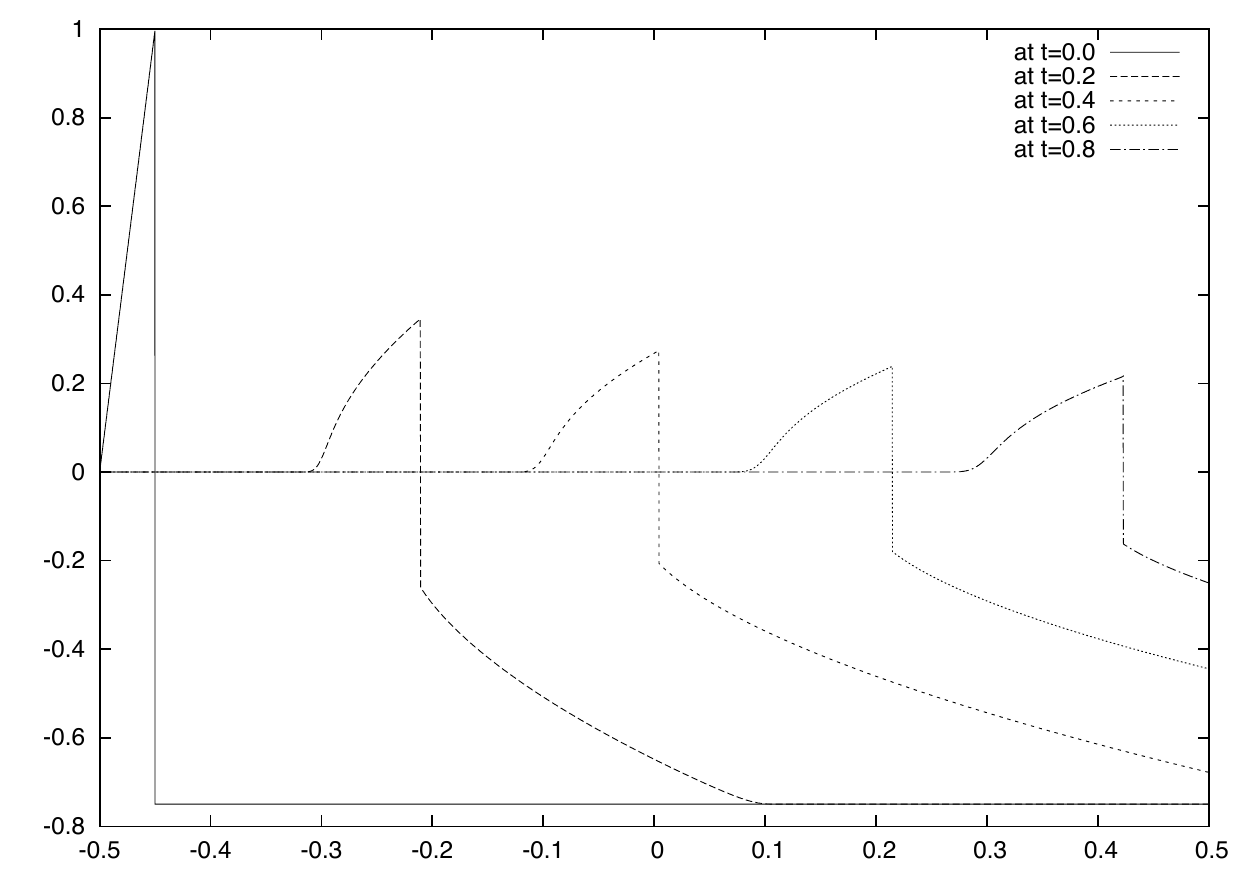}
\includegraphics[width=8.0cm]{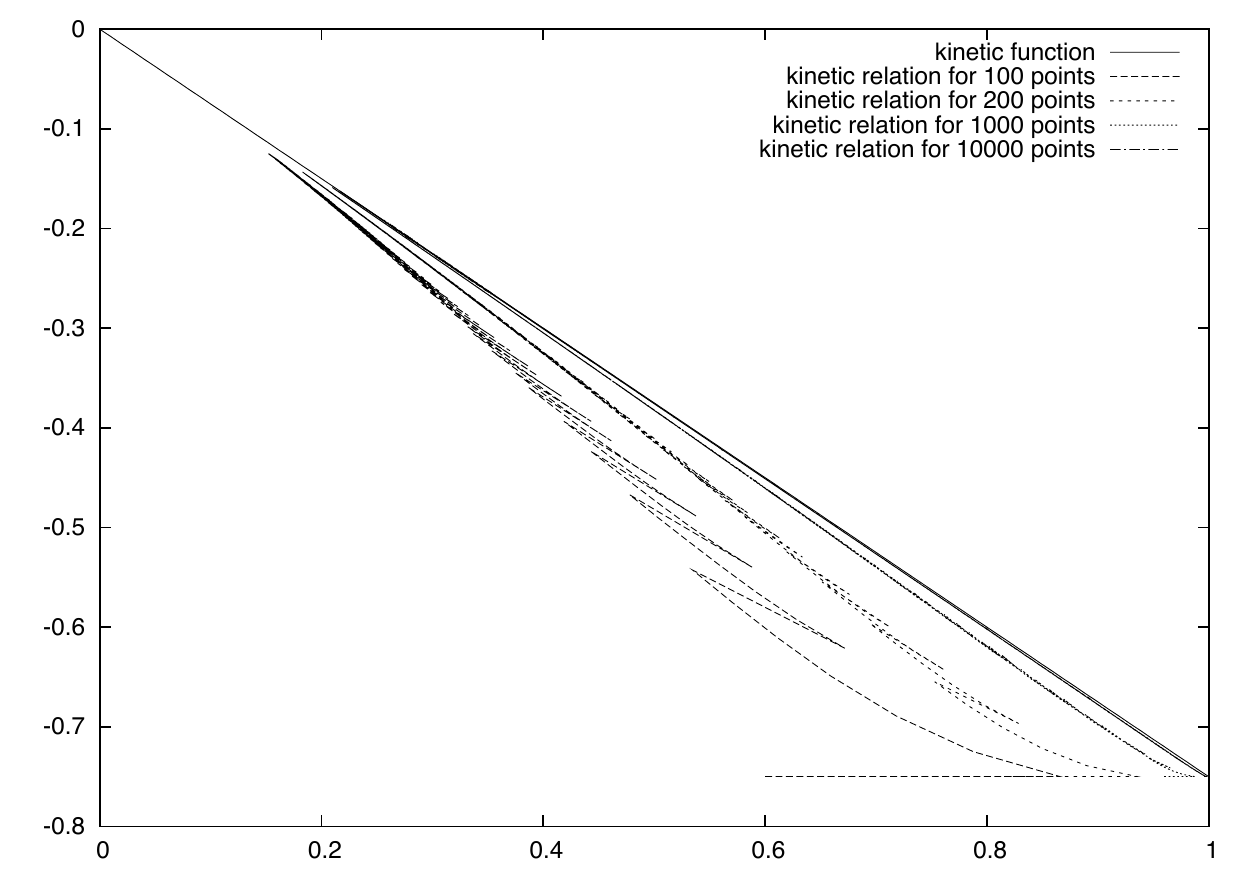}
\vspace{-0.45cm}
\caption{Test F - Numerical kinetic relation} 
\label{Test6}
\end{center}
\end{figure}

{\bf Test F.} We now study how the kinetic relation
$u_R=\varphi^\flat(u_L)$ is computed. On Figure \ref{Test6} (right-hand figure), we plot points whose
horizontal coordinates (respectively vertical coordinates) correspond to the left-hand 
(resp. right-hand) traces around the reconstructed cell. The initial data 
allows us to cover a large range of value: 
$$
u_0(x)=
\left\{
\begin{array}{rl}
0, & x<0.5,\\
1 + 20(x+0.45), & 0.5<x< 0.45,\\
-0.75, & x>-0.45.
\end{array}
\right.
$$
The left-hand figure represents the solution at different times with
$\Delta x=0.0002$.\\

We clearly observe the convergence of the numerical kinetic relation
towards the prescribed one. This a strong test to validate the proposed method. 


{\bf Test G.} In the course of designing the scheme proposed in the previous
section we tried several variants. We report here one such scheme that
is very similar to the proposed scheme, but which {\sl does not}
converge to exact nonclassical solutions. This is due to the fact that
small oscillations are generated in the scheme which are in competition with
the dissipation mechanisms described by the prescribed kinetic
function. 

\begin{figure}[htbp]
\begin{center}
\includegraphics[width=9.0cm]{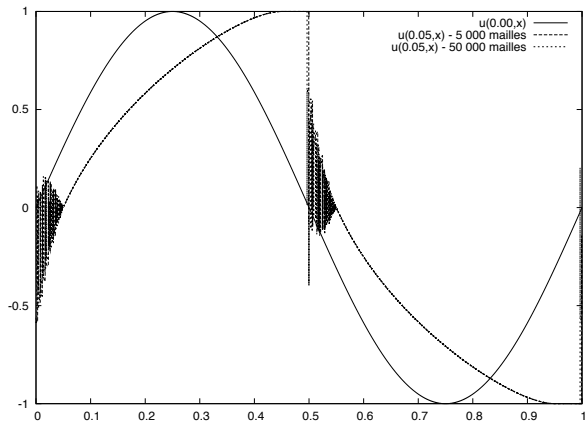}
\caption{Another version of the scheme} 
\label{BadScheme}
\end{center}
\end{figure}

The variant is designed for the concave-convex flux $f(u) = u^3 +
u$. The only difference with the scheme developed above is that it
performs the reconstruction in $\mathcal{C}_j$ with $u_{j,l}^n =
u_{j-1}^n$ (instead of $\varphi^{-\flat}(u_{j+1}^n$) and $u_{j,r} =
\varphi^\flat (u_{j-1}^n)$. This is equivalent in the case of a {\em
pure} nonclassical shock (Test B) but different in the general case. \\ 
Figure \ref{BadScheme} presents the solution obtained for the same
initial value as in Test E. Oscillations are generated because the
reconstruction is not constrained enough in this version of the
scheme.\\


\section{Concluding remarks} 

In this paper we have introduced a new numerical strategy for
computing nonclassical solutions to nonlinear hyperbolic conservation
laws. The method is based on a reconstruction technique performed in
each computational cell which may exhibit a nonclassical
shock. Importantly, the whole algorithm is {\sl conservative} and
propagates any admissible nonclassical discontinuity exactly. The
convergence of the proposed method was demonstrated numerically for
several test-cases. This new approach brings a new perspective on the
numerical approximation of nonclassical shocks and kinetic
functions. The efficiency of the method is clearly demonstrated in the
present paper, and we refer to the follow-up paper \cite{BCLL-two} for
various extensions and applications. Among the questions of interest
we can mention the total variation bounds and the hyperbolic systems of conservation laws, 
the application to real materials undergoing phase transitions, as well
as the extension to higher-order schemes. 


\end{document}